\documentclass[review,onefignum,onetabnum]{siamonline220329}

\usepackage{bm}
\usepackage{amsmath}
\usepackage{amssymb}
\usepackage{subfigure}
\usepackage{overcite}
\usepackage{footnpag}			      	
\usepackage{booktabs}
\usepackage{algorithmic}
\usepackage[hang,flushmargin]{footmisc}

\newtheorem*{remark}{Remark}

\DeclareMathOperator*{\argmin}{arg\,min}

\begin{document}

\title{Rapid GPU-Assisted Search and Parameterization-Based Refinement and Continuation of Connections between Tori in Periodically Perturbed Planar Circular Restricted 3-Body Problems\thanks{Submitted to the editors DATE. Part of this article appeared as manuscript no. AAS 21-349 in the non peer-reviewed proceedings of the 31st AAS/AIAA Space Flight Mechanics Meeting, Virtual, February 1-3, 2021.
\funding{The first author was supported by a NASA Space Technology Research Fellowship under grant no. 80NSSC18K1143. Part of the writing of this article was supported by the NSF under award no. DMS-2202994}}}

\author{Bhanu Kumar\thanks{School of Mathematics, Georgia Institute of Technology, 686 Cherry Street NW, Atlanta, GA 30332 (\email{bkumar30@gatech.edu}, \email{rafael.delallave@math.gatech.edu}).}  
\and Rodney L. Anderson\thanks{Jet Propulsion Laboratory, California Institute of Technology, 4800 Oak Grove Drive, Pasadena, CA 91109 (\email{Rodney.L.Anderson@jpl.nasa.gov}).}
\and Rafael de la Llave\footnotemark[2]}

\headers{Rapid GPU and Parameterization-Based Search of Connections between Tori}{B. Kumar, R.L. Anderson, and R. de la Llave}

\maketitle

\begin{abstract} 
When the planar circular restricted 3-body problem (PCRTBP) is periodically perturbed, as occurs in many useful astrodynamics models, most unstable periodic orbits persist as whiskered tori. Intersections between stable and unstable manifolds of such tori provide natural heteroclinic pathways enabling spacecraft to greatly modify their orbits without using propellant. However, the 2D Poincar\'e sections used in PCRTBP studies no longer work to find these intersections. Thus, in this study, we develop new fast methods to search for and compute such heteroclinics. First, the dynamics are used to restrict the intersection search to only certain manifold subsets, greatly reducing the required computational effort. Next, we present a massively parallel procedure for carrying out this search by representing the manifolds as discrete meshes and adapting methods from computer graphics collision detection algorithms. Implementing the method in Julia and OpenCL, we obtain a 5-7x speedup by leveraging GPUs versus CPU-only execution. Finally, we show how to use manifold parameterizations to refine the approximate intersections found in the mesh search to very high accuracy, as well as to numerically continue the connections through families of tori; the families' Whitney differentiability enables interpolation of needed parameterizations. The ability to very rapidly find a heteroclinic intersection between tori of fixed frequencies thus allows the systematic exploration of intersections for tori of nearby frequencies as well, yielding a variety of potential zero-fuel spacecraft trajectories. We demonstrate the tools on the Jupiter-Europa planar elliptic RTBP.

\end{abstract}



\section{Introduction}

Numerous prior studies have used the stable and unstable manifolds of unstable planar periodic orbits, both around libration points as well as at resonances, as an efficient tool for multi-body mission design in the planar circular restricted 3-body problem (PCRTBP). For instance, Anderson and Lo\cite{Anderson2010,Anderson2011} studied intersections of the manifolds of resonant periodic orbits in the Jupiter-Europa system as a mechanism of resonance transitions. The book of Koon et al.\cite{KoLoMaRo} describes the Poincar\'e section method of finding intersections of manifolds between L1 and L2 libration point planar Lyapunov orbits, and shows how to use the resulting regions to construct trajectories with arbitrary itineraries between different realms of the PCRTBP model. As the PCRTBP phase space is 4-dimensional, fixing an energy level restricts the dynamics to a 3D submanifold, and the Poincar\'e section further reduces the dimensionality of the system to 2D. Since the manifolds of the periodic orbits are 2D cylinders in the full phase space, taking the Poincar\'e section reduces the manifolds to 1D curves in the section, so the problem of finding connections between periodic orbits for space missions to follow reduces to finding the intersection of two 1D manifold curves in a 2D plane.

While this method of intersecting 1D curves is useful for the 4-dimensional PCRTBP, it is not applicable to systems with higher dimensional phase spaces. For instance, in the 6D phase space of the spatial CRTBP, fixing an energy level and taking a Poincar\'e section results in a 5D energy level and 4D section. In the case of a time-varying periodic forcing of the PCRTBP -- as occurs when an effect neglected in the PCRTBP (e.g. a third large body) is reintroduced into the model to improve accuracy  -- the phase space becomes 5-dimensional, considering time as a state variable. As energy is no longer conserved in non-autonomous systems, taking a Poincar\'e section again leaves us with a 4D space to be explored. Furthermore, in these higher dimensional systems, unstable periodic orbits are no longer the main dynamical structures of interest for zero-fuel orbit transfers;  2D manifolds of periodic orbits do not generically intersect each other in a 5D phase space or energy level. Instead, as predicted by Arnol'd\cite{arnold1963}, unstable quasi-periodic orbits (also known as \emph{whiskered tori}) and their manifolds are the most common objects of sufficient dimensionality such that one can expect generic intersections. 

In prior work\cite{kumar2022}, we described how most PCRTBP unstable periodic orbits persist as 2D whiskered tori in the 5D phase space of periodically-perturbed PCRTBP models, at least when the strength of the perturbation is sufficiently small; we also developed and implemented very efficient and precise methods of computing such tori and their invariant manifolds. By considering stroboscopic maps instead of the continuous-time flow, we can reduce the dimensionality of the system by 1 so that the 2D quasi-periodic orbits become invariant 1D tori (circles) in the 4D stroboscopic map phase space $(x,y,p_{x},p_{y})$; these invariant circles have 2D cylindrical stable and unstable manifolds as the PCRTBP unstable periodic orbits did. However, due to the absence of an energy integral, manifold intersections in the perturbed system will occur at isolated points, rather than along continuous trajectory curves. Hence, a different method of computing homoclinic and heteroclinic connections in the map's 4D phase space is required. 

The goal of this study is to develop new and computationally fast methods and tools for computing intersections between the previously described quasi-periodic solutions' manifolds. Such intersections provide trajectories in higher-accuracy, periodically-perturbed PCRTBP models that a spacecraft can follow to traverse the system phase space without using any propellant. Our algorithms are designed in large part for implementation on modern graphics processing units (GPUs) to greatly accelerate their performance; a speedup of 5-7x is achieved in the case of our implementation. Since this paper blends together concepts from mathematics (KAM theory, parameterization methods for invariant manifolds, Whitney differentiability), computer graphics (collision detection algorithms), and GPU computing, we have included short descriptions of these tools, which the reader can skip depending on background. 

In this paper, we start with a brief overview of GPU computing capabilities and paradigms as well as some background on collision detection methods from computer graphics, both of which are used in this work. Then, after defining some dynamical models and maps, we give a summary of our previously developed\cite{kumar2022} parameterization method for the computation of whiskered tori and their manifolds in periodically-perturbed PCRTBP models. Next, we describe the method of layers for restricting the homoclinic and heteroclinic connection search to appropriate subsets of the two manifolds of interest. Once these subsets are identified, we can computationally represent each manifold as a 2D mesh of points in the 4D stroboscopic map phase space. With these meshes representing the manifolds, we then develop a heavily parallel algorithm for detecting and computing intersections of these meshes in 4D space; we then describe how to implement this method using the Julia programming language and OpenCL, taking advantage of the capabilities of modern GPUs to greatly speed up the algorithm execution time.  Finally, we take the approximate intersections of manifolds computed from the mesh intersection search, and show how to use our manifold parameterizations to refine the intersections to high precision and then continue them through 1-parameter families of tori of varying frequencies. We demonstrate the use of our methods by applying them to the search for heteroclinic connections between resonances in the Jupiter-Europa planar elliptic RTBP.

\section{Background}
\subsection{An Overview of GPU Computing}

Graphics processing units (GPUs) are special-purpose computer processors originally designed for executing 3D graphics-related computations\cite{gpucomputing}. Nevertheless, as GPU capabilities and availability have grown, it has become possible to use them for many other computational tasks as well. GPUs excel at tasks with large and highly parallel computational requirements, as the GPU processes blocks of many elements in parallel using the same program. While a single-program multiple-data (SPMD) programming model is supported on GPUs, due to the lock-step execution of the program on multiple data elements, it is necessary to evaluate both sides of any code branches for all elements in a block of data. Hence, GPUs are best suited to straight-line programs which are mostly written in a single-instruction, multiple data (SIMD) style. Flow control should be kept to a minimum. Since not all algorithms run well on GPUs, it is not always true that a problem can be implemented on them in a performant manner; even when a GPU-based solution does exist, the required programs need to be carefully designed. 

On a practical level, there are two main toolkits used for programming GPUs. The most common is Nvidia CUDA\cite{cuda}, which works only with Nvidia GPUs. The other is OpenCL,\cite{khronos} which is an open standard that provides a programming language and APIs for a variety of SIMD-capable devices, including both AMD and Nvidia GPUs as well as CPUs. OpenCL implementations exist for many different platforms, including Windows, MacOS, and Linux x64. The basic programming concepts of CUDA and OpenCL are very similar; in this study we used OpenCL, so we now briefly review the OpenCL model and concepts used in this work. 

In OpenCL, the fundamental task is the programming of kernels. A \emph{kernel} is a program which executes on the OpenCL device, such as a GPU. When a kernel is submitted for execution, a 1D, 2D, or 3D space of indices is defined, called an \emph{NDRange} (we only use a 1D NDRange in this study); each index corresponds to a separate execution of the kernel. These kernel instances are referred to as \emph{work items}, with each work item identified by a unique global ID (a nonnegative integer if NDRange is 1D) corresponding to an index in NDRange. Each work item executes the same kernel, but can access different data in memory by using its unique index value. Work items are grouped into \emph{work groups}; all work items in a work group execute concurrently. Finally, the work groups taken together form the NDRange. The key advantage of GPUs is the massive number of threads they have, usually on the order of thousands, which allows them to execute large numbers of work items in parallel. 

The OpenCL memory model has various categories of memory stored on the device, accessible to different parts of the program. \emph{Private} memory is accessible only to a single work item, and \emph{local} memory is accessible to all work items in the same work group. We do not use local memory in the programs written for this study. Finally, \emph{global} memory is accessible to all work items. The use of global memory can be an issue if two work items try to access the same memory location at the same time; to solve this problem, \emph{atomic} functions are useful. These functions receive a pointer to a 32 bit integer or floating point number stored in global or local memory, modify the value there, and return the old value. The key is that atomic functions execute so that when one work item is carrying out an atomic operation at a pointer location, the other threads must wait for that work item to complete the operation. Some such functions are atomic add, subtract, increment, min, and max. 

\begin{figure}
\begin{centering}
\includegraphics[width=0.99\columnwidth]{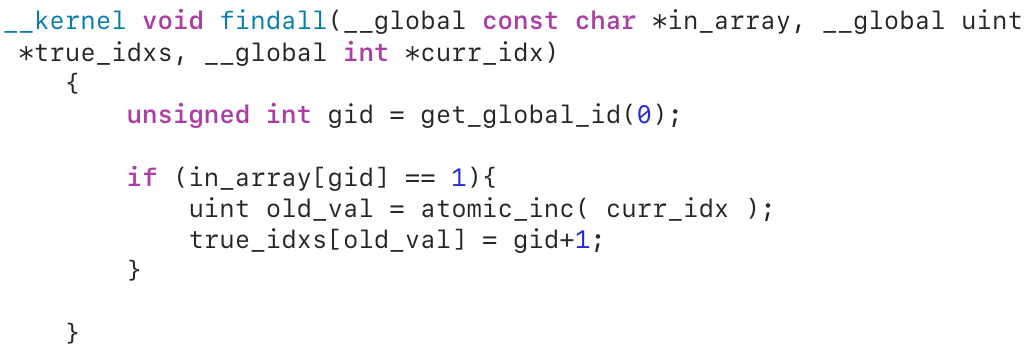}
\caption{ \label{fig:findall} Example of OpenCL findall kernel}
\end{centering}
\end{figure}

An example of a simple OpenCL kernel is given in Figure \ref{fig:findall}. This program takes as input an array of zeros and ones (denoted as in\_array) and finds the in\_array entries which have a value of 1, saving their indices in the array true\_idxs. It also takes a pointer curr\_idx to an integer which is initialized to 0 before the kernel is run. The first line after the kernel declaration stores the work item global ID as gid. Then, the kernel checks if the gid entry of in\_array is 1, so that each work item checks a different entry of in\_array. If the gid entry is true, an atomic increment is applied to curr\_idx, with the pre-increment value of curr\_idx stored in old\_val; each work item stores its value of old\_val in private memory, separately from all other work items. curr\_idx is in global memory, so it is shared between all work items. Hence, the first work item to apply atomic\_inc sets its private old\_val to 0 and the shared curr\_idx to 1, the next work item to apply atomic\_inc sets its old\_val to 1 and curr\_idx to 2, and so on. This way, each work item with in\_array[gid] equal to 1 gets a unique consecutive value of old\_idx. Finally, the gid values with in\_array[gid] equal to 1 are stored in the true\_idxs array at the indices old\_idx by the corresponding work items (we add 1 to gid since we later use the true\_idxs array in Julia, which has one-based indexing).  

\subsection{Collision Detection} \label{colDetectSection}

Given two sets $A$ and $B$, each containing some finite number of geometric objects, the problem of collision detection simply is that of determining whether any object from $A$ intersects any object(s) from $B$, and if so, which pairs of objects intersect and where they do so. Computational methods for collision detection are used in a variety of applications, including video game development, animation, robotics, computer aided design, and physics simulations \cite{collisionDetection, gpuGems3}. Indeed, as will be shown in Section \ref{gpuSearchSection}, such methods are also applicable to the problem of finding heteroclinic connections in dynamical systems. 

If one were to solve the collision detection problem by checking each object in $A$ against each object in $B$ for an exact determination of whether/where they intersect, this would require $|A| |B|$ tests; even if $A=B$ as is often the case in computer graphics (though not in this paper), there would still be $\frac{|A|(|A|-1)}{2}$ pairwise intersection tests required. Either way, the cost of such a collision detection algorithm is quadratic, which can quickly become infeasible if $|A|$ and $|B|$ are large (as will be the case in this paper) and the pairwise intersection test itself is not extremely computationally cheap. Unless the objects of $A$ and $B$ are all very simple, though (e.g. spheres), the exact pairwise intersection test is generally expensive. 

It is hence necessary to find a way to reduce the number of exact pairwise tests done. To this end, most collision detection algorithms are comprised of two phases, called broad phase and narrow phase.\cite{gpuGems3} Since in practice most pairs of objects do not intersect, the broad phase uses very cheap computations to eliminate most of these non-intersecting pairs from consideration. Usually this involves dividing objects into groups such that only objects within the same group can intersect, and/or putting pairs of objects through a very simple test which, if failed, verifies non-intersection. The broad phase outputs a set of potentially colliding pairs of objects, to which the narrow phase then applies the more costly exact intersection test. We will apply these ideas to find collisions between pieces of invariant manifolds in Section \ref{gpuSearchSection}. 

\section{Dynamical Models and Equations of Motion} \label{modelsection}
\subsection{Planar Circular Restricted 3-body Problem}
The planar circular restricted 3-body problem (PCRTBP) models the motion of a spacecraft under the gravitational influence of two large bodies which revolve in a circular Keplerian orbit around their barycenter. The large bodies are called the primary body of mass $m_{1}$ and secondary body of mass $m_{2}$ (together referred to as the primaries), with mass ratio $\mu = \frac{m_{2}}{m_{1}+ m_{2}}$. After normalizing units, the distance between $m_{1}$ and $m_{2}$ becomes 1, $\mathcal{G}(m_{1}+m_{2})=1$, and their period of revolution becomes $2 \pi$. Using the usual synodic coordinate system with the primaries on the $x$-axis and origin at their barycenter, the equations of motion are Hamiltonian of form \cite{celletti}
\begin{equation} \label{pcrtbpH_EOM} \dot x = \frac{\partial H_{0}}{\partial p_{x}} \quad \dot y = \frac{\partial H_{0}}{\partial p_{y}} \quad \quad \dot p_{x} = -\frac{\partial H_{0}}{\partial x} \quad \dot p_{y} = -\frac{\partial H_{0}}{\partial y} \end{equation}
\begin{equation}  \label{pcrtbpH} H_{0}(x,y,p_x,p_{y})= \frac{p_{x}^{2}+p_{y}^{2}}{2} + p_{x}y -p_{y}x - \frac{1-\mu}{r_{1}} - \frac{\mu}{r_{2}} \end{equation}
where $r_{1} = \sqrt{(x+\mu)^{2} + y^{2}}$ and $r_{2} = \sqrt{(x-1+\mu)^{2} + y^{2}} $ are the distances from the spacecraft to $m_{1}$ and $m_{2}$, respectively. We assume that the spacecraft moves in the same plane as the primaries. The Hamiltonian in Equation \eqref{pcrtbpH} is autonomous and is thus an integral of motion.

\subsection{Periodic Perturbations of the PCRTBP}\label{perturbedpcrtbpsection}
Many effects not modeled in the PCRTBP act as a time-periodic forcing on the spacecraft; these can be represented as periodic perturbations of the PCRTBP. In some cases these perturbations can have a significant effect on the trajectory and thus should be accounted for even during initial mission design. Many such perturbations are also Hamiltonian; in this case, the equations of motion become 
\begin{equation} \label{H_EOM} \dot x = \frac{\partial H_{\varepsilon}}{\partial p_{x}} \quad \dot y = \frac{\partial H_{\varepsilon}}{\partial p_{y}} \quad \quad \dot p_{x} = -\frac{\partial H_{\varepsilon}}{\partial x} \quad \dot p_{y} = -\frac{\partial H_{\varepsilon}}{\partial y} \quad \quad \dot \theta_{p} = \Omega_{p} \end{equation}
\begin{equation}  \label{perturbed_H} H_{\varepsilon}(x,y,p_x,p_{y}, \theta_{p})= H_{0}(x,y,p_x,p_{y})+ H_{1}(x,y,p_x,p_{y}, \theta_{p}; \varepsilon)\end{equation}
where $\theta_{p} \in \mathbb{T}$ is an angle considered modulo 2$\pi$, $H_{0}$ is the PCRTBP Hamiltonian from Equation \eqref{pcrtbpH},  $H_{1}$ is the time-periodic perturbation satisfying $H_{1}(x,y,p_x,p_{y}, \theta_{p}; 0)=0$, and $\varepsilon > 0$ and $\Omega_{p}$ are the perturbation parameter and perturbation frequency, respectively. $\varepsilon$ signifies the strength of the perturbation, with $\varepsilon=0$ being the unperturbed PCRTBP, and $\Omega_{p}$ is a known constant frequency. The forcing from $H_{1}$ is $2\pi / \Omega_{p}$ periodic, with $\theta_{p}$ being the perturbation phase angle. The perturbed system's phase space is thus $(x,y,p_x,p_{y}, \theta_{p}) \in \mathbb{R}^{4} \times \mathbb{T}$. In general, the Hamiltonian in Equation \eqref{perturbed_H} is not an integral of motion when $\varepsilon \neq 0$. 

\subsection{Planar Elliptic Restricted 3-body Problem}
There are many periodically perturbed PCRTBP models that are of interest, such as the bicircular problem \cite{simo1995bicircular}, the quasi-bicircular problem \cite{andreu1998quasi}, and the Hill restricted 4-body problem \cite{scheeres1998restricted}. Another well-known periodically perturbed PCRTBP model is the planar elliptic restricted 3-body problem (PERTBP), which we use in this study for numerical demonstration of our tools.  In the PERTBP, $m_{1}$ and $m_{2}$ revolve around their barycenter in an elliptical Keplerian orbit of nonzero eccentricity $\varepsilon > 0$. All other assumptions are the same as the PCRTBP. The length unit is defined such that the $m_{1}$-$m_{2}$ orbit semi-major axis is 1, and the period of the primaries' orbit remains $2\pi$. This implies that the perturbation frequency $\Omega_{p} = 1$; hence, we can take $\theta_{p} = t$ modulo $2\pi$. 

We use the same PERTBP model used by Hiday-Johnston and Howell \cite{hiday1994transfers}, except for a coordinate change from position-velocity to position-momentum coordinates and a restriction of the dynamics to the $xy$-plane. Again, the coordinate system is defined with the primaries on the $x$-axis and origin at their barycenter. However, the distance from $m_{1}$ to $m_{2}$ is now time-periodic; we let $t=0, 2\pi, \dots$ be the times of periapse of their orbit. This is different from the well-known pulsating coordinates of Szebehely \cite{szebehely1969}. The equations of motion are given by Equation \eqref{H_EOM} with time-periodic Hamiltonian 
\begin{equation}  \label{pertbpH} H_{\varepsilon}(x,y,p_x,p_{y},t)= \frac{p_{x}^{2}+p_{y}^{2}}{2} + n(t)(p_{x}y -p_{y}x) - \frac{1-\mu}{r_{1}} - \frac{\mu}{r_{2}} \end{equation}
where $r_{1} = \sqrt{(x+\mu(1-\varepsilon \cos E(t)))^{2} + y^{2}} $ and $ r_{2} = \sqrt{(x-(1-\mu)(1-\varepsilon \cos E(t)))^{2} + y^{2}}$. $E(t)$ is the $2\pi$-periodic eccentric anomaly of the elliptical $m_{1}$-$m_{2}$ orbit, and can be computed by solving the usual Kepler's equation $M = E - \varepsilon \sin E$ as described in Bate, Mueller, and White \cite{bmw} ($M=t$ in our case). $n(t)$ is the time derivative of the $m_{1}$-$m_{2}$ true anomaly. From Equations \eqref{H_EOM} and \eqref{pertbpH}, the momenta are related to velocity by $p_{x}=\dot x-n(t)y$ and $p_{y}=\dot y+n(t)x$. We will use the Jupiter-Europa PERTBP in this study, where $\mu \approx 2.527 \times 10^{-5}$ and $\varepsilon = 0.0094$.

\section{Summary of Previous Results on Computing Tori and Manifolds} \label{backgroundsection}

In this section, we briefly review the results of the parameterization methods developed in Kumar et al.\cite{kumar2022}, which themselves were inspired by methods described in Haro et al. \cite{haroetal} and Zhang and de la Llave\cite{zhang}. 

\subsection{Stroboscopic Maps} \label{stroboscopic}

The quasi-periodic orbits of interest in periodically-perturbed PCRTBP models lie on 2D invariant tori in the 5D extended phase space $(x,y,p_{x},p_{y}, \theta_{p})$. These invariant tori can be parameterized as the image of a function of two angles $K_2 : \mathbb{T}^{2} \rightarrow \mathbb{R}^{4} \times \mathbb{T}$. Any quasi-periodic trajectory $\bold{x}(t)$ lying on this torus can be expressed as
	\begin{equation} \bold{x}(t) = K_2(\theta, \theta_p)  \quad \quad \quad \theta = \theta_{0}+\Omega_1 t, \quad  \theta_p = \theta_{p,0}+\Omega_p t \end{equation}
where the initial condition $\bold{x}(0)$ determines $\theta_{0}$ and $\theta_{p,0}$. $\theta_{p}$ and $\Omega_p$ are the perturbation phase and frequency defined earlier, respectively. Defining the stroboscopic map $F: \mathbb{R}^{4} \times \mathbb{T} \rightarrow \mathbb{R}^{4} \times \mathbb{T}$ as the time-$2\pi/\Omega_p$ mapping of extended phase space points by the equations of motion, we have
	 \begin{equation} \label{invariance_fK2} F(K_2(\theta, \theta_p)) = K_2(\theta+\omega, \theta_p), \text{ where  } \omega = 2\pi \Omega_1/\Omega_p \end{equation} 
	 since the angle $\theta_p$ increases by $2\pi$ in the time $2\pi/\Omega_p$. Since the value of $\theta_{p}$ does not change under the map $F$, one can fix a value of $\theta_p$ and define $K(\theta) = K_2(\theta, \theta_p)$ (without loss of generality, we choose $\theta_{p}=0$ in this study). Then, Equation \eqref{invariance_fK2} becomes 
\begin{equation} \label{invariance} F(K(\theta)) = K(\theta+\omega)  \end{equation} 
Ignoring the invariant $\theta_{p}$ component of the extended phase space and making a slight abuse of notation, we have $F:\mathbb{R}^{4} \rightarrow \mathbb{R}^{4}$ and $K:\mathbb{T} \rightarrow \mathbb{R}^{4}$. Equation \eqref{invariance} implies that $K$ is an invariant 1D torus of $F$. Hence, basing our study on the stroboscopic map $F$ is more efficient than solving for tori invariant under the flow of the ODE, since we reduce the phase space dimension from 5D to 4D and the dimension of the unknown invariant tori from 2D to 1D. This greatly reduces the computational cost of representing functions and objects, as well as of calculations involving them. Thus, we use this approach in the remainder of this study.  

\subsection{Parameterization Methods for Tori, Bundles, and Manifolds} \label{quasiNewton}

With the stroboscopic map $F$ defined, our first goal is to find solutions $K(\theta)$ of Equation \eqref{invariance} parameterizing the tori. The rotation number $\omega=2\pi \Omega_1/\Omega_p$ is generally known; for instance, this is the case whenever the invariant torus being solved for comes from a known PCRTBP periodic orbit. 

The quasi-Newton method developed in our previous work\cite{kumar2022} for solving Equation \eqref{invariance} adds an extra equation to be solved.  In particular, in addition to solving for  $K(\theta)$, we simultaneously solve for matrix-valued periodic functions $P(\theta)$, $\Lambda(\theta): \mathbb{T} \rightarrow \mathbb{R}^{4 \times 4}$ satisfying
		 \begin{equation}  \label{bundleEquations} DF(K(\theta)) P(\theta) = P(\theta+\omega) \Lambda(\theta) \end{equation} 
$P(\theta)$ and $\Lambda(\theta)$ are the matrices of bundles and of Floquet stability, respectively; for each $\theta \in \mathbb{T}$, the columns of $P(\theta)$ are comprised of the tangent, symplectic conjugate center, stable, and unstable directions of the torus at the point $K(\theta)$, in that order, while $\Lambda$ has the form
\begin{equation}  \label{bundles_equation}
\Lambda(\theta)=\begin{bmatrix}
1 &  T(\theta)   & 0 & 0 \\ 0 &  1   & 0 & 0 \\ 0 & 0  & \lambda_s  & 0 \\ 0 &  0 & 0 & \lambda_u \end{bmatrix}
 \end{equation}
where $T:\mathbb{T} \rightarrow \mathbb{R}$ and $\lambda_s, \lambda_u \in \mathbb{R}$ are constants with $\lambda_s<1$ and $\lambda_u > 1$. 

As it turns out, solving simultaneously for $K$, $P$, and $\Lambda$  not only gives more information than solving for $K$ alone, but is also more efficient.  In fact, discretizing the torus on a grid of $N$ points, our quasi-Newton method requires $O(N)$ storage and an operation count of only $O(N \log N)$, as compared to $O(N^{3})$ operations for $K$-only methods. Moreover, the algorithm is expressed in terms of a few vector operations, each of which can be implemented easily in a high level language (such as Julia) either in Fourier space or in real space, with FFTs allowing for easy and fast conversions between these two representations. We refer the reader to Kumar et al\cite{kumar2022} for full details. Given a sufficiently accurate initial guess, each step of the quasi-Newton method reduces the error in Equations \eqref{invariance}-\eqref{bundleEquations} to roughly the square of the original error (as in the Newton method); it can thus be used for numerical continuation of tori from the PCRTBP ($\varepsilon=0$) to perturbed models with $\varepsilon>0$. Furthermore the algorithm is backed by \emph{a posteriori} theorems \cite{FontichLS09a} which show that close to numerical results with small residual, there are true solutions. Besides its theoretical interest, the a posteriori theorems specify the condition numbers one needs to monitor to be confident of the results.

Once the $F$-invariant tori and their center, stable, and unstable torus bundles are computed, we next wish to accurately compute stable and unstable manifolds of the tori under $F$. As the invariant tori are 1D, their stable and unstable manifolds will be 2D and topologically equivalent to $\mathbb{T} \times \mathbb{R}$. We can hence parameterize the manifolds as the image of a function $W(\theta, s)$ satisfying 
\begin{equation}  \label{invariancequationfinal}   F(W(\theta, s)) = W(\theta + \omega, \lambda s) \end{equation}
where $(\theta, s) \in \mathbb{T} \times \mathbb{R}$, and $\lambda$ is the stable $\lambda_{s}$ or unstable $\lambda_{u}$ multiplier from $\Lambda$, depending on which manifold we are solving for. To solve Equation \eqref{invariancequationfinal}, we express $W$ as a Fourier-Taylor series
\begin{equation}  \label{series} W(\theta, s) = \sum_{k \geq 0} W_{k}(\theta)s^{k} = K(\theta) + \sum_{k \geq 1} W_{k}(\theta)s^{k}  \end{equation}
where $K(\theta)$ is the invariant circle whose manifold we are trying to compute. The $s^{0}$ term of $W$ is $W_{0}(\theta)=K(\theta)$, and the linear term $W_{1}(\theta)$ is the stable or unstable bundle known from the third or fourth column of $P$. The higher order $W_{k}(\theta)$ terms can then be solved for recursively, as described in our previous work\cite{kumar2022}. Note that $s=0$ corresponds to the base invariant torus.

Note that Equations \eqref{invariance}-\eqref{bundleEquations}, as well as Equation \eqref{invariancequationfinal} are underdetermined. Changing phase (replacing $\theta $ with $ \theta + \rho$, $\rho \in \mathbb{T}$) and/or scales (replacing $s $ with $ Ls$,  $L \in \mathbb{R}$) in $K$, $P$, $\Lambda$, and $W$ also leads to solutions of those equations. This corresponds to changing variables in the manifold used as the domain of the parameterizations, and does not affect the tori or invariant manifolds computed. Both of those underdeterminacies can be eliminated by fixing normalizations. The normalized solutions of \eqref{invariance}-\eqref{bundleEquations} and \eqref{invariancequationfinal} are locally unique, so we can compare the results of different implementations of the algorithm. One can also discuss smoothness of the normalized solutions with respect to parameters. From the numerical point of view, the scale normalization affects the round-off error of the algorithms while the normalization of the phase is irrelevant. See Kumar et al \cite{kumar2022} and Section \ref{W1fSection} for more details.

\subsection{Using Parameterizations for Manifold Globalization} 

We have found that the torus stable and unstable manifold parameterizations $W(\theta, s)$ described in the previous section approximate the manifolds very accurately, but only within some finite range of $s$ values. To find this range of values, we first choose an error tolerance, say $E_{tol} = 10^{-5}$ or $10^{-6}$. We then  find the largest $D \in \mathbb{R}^{+}$ such that for all $s$ satisfying $|s| \leq D $, 
\begin{equation} \max_{\theta \in \mathbb{T}} \|F(W(\theta, s)) - W(\theta+\omega, \lambda s)\| < E_{tol} \end{equation}
We refer to the set $\mathbb{T} \times (-D,D)$ as the fundamental domain of $W(\theta, s)$. Once $D$ is computed, we wish to extend $W(\theta, s)$ and compute manifold points for $s$ values outside the fundamental domain. This step is called globalization. Repeatedly applying Equation \eqref{invariancequationfinal}, we have that $F^{k}(W(\theta,s)) = W(\theta+k \omega, \lambda^{k}s)$, where superscript $k$ on $F$ refers to function composition. From this, we can derive two relations:
\begin{gather} \label{globou} W(\theta, s) = F^{k}(W(\theta-k \omega, \lambda^{-k}s)) \\
\label{globos} W(\theta, s) = F^{-k}(W(\theta+k \omega, \lambda^{k}s)) \end{gather}
These equations can be used to evaluate $W(\theta, s)$ for $s$ values outside the fundamental domain. If $W$ is an unstable manifold with $|\lambda| > 1$, then take $k \geq 0$ such that $| \lambda^{-k}s | < D$ and evaluate Equation \eqref{globou}, using the Fourier-Taylor series for $W$ to compute $W(\theta-k \omega, \lambda^{-k}s)$ . Similarly, if $W$ is a stable manifold with $|\lambda| < 1$, then take $k \geq 0$ such that $| \lambda^{k}s | < D$ and evaluate Equation \eqref{globos}. The mapping by $F^{k}$ or $F^{-k}$ is just computed using numerical integration. 

As a final remark, note that the manifold parameterization along with Equations \eqref{globou}-\eqref{globos} give very detailed information not just on the manifold itself, but also on its tangent directions at each point. These tangents can be found by simply differentiating $W(\theta,s)$ with respect to $\theta$ and $s$, which knowledge of the Fourier-Taylor series for $W$ makes possible. As we show in Section \ref{refinementSection}, this is extremely useful when refining approximate manifold heteroclinic intersections, as it enables the use of Newton methods and implicit function theorems. 

\section{The Method of Layers for Restricting the Connection Search} \label{layerSection}

After computing the stable and unstable manifolds of unstable invariant tori (invariant circles) in a periodically-perturbed PCRTBP, a natural next step is to search for heteroclinic connections between them. Henceforth, let $W^{u}_{1}(\theta_{u}, s_{u})$ and $W^{s}_{2}(\theta_{s}, s_{s})$ represent the unstable and stable manifolds of stroboscopic map invariant circles 1 and 2, respectively. Heteroclinic connections from circle 1 to circle 2 occur when the images of $W^{u}_{1}$ and $W^{s}_{2}$ intersect in $(x,y,p_{x}, p_{y})$ space. This means we need to find $(\theta_{u}, s_{u})$ and $(\theta_{s}, s_{s})$ such that 
\begin{equation} \label{intersectcondition} W^{u}_{1}(\theta_{u}, s_u) = W^{s}_{2}(\theta_{s}, s_{s}) \end{equation}
In order to solve Equation \eqref{intersectcondition}, it would help to be able to restrict our solution search to only certain regions of the $(\theta_{u}, s_{u}, \theta_{s}, s_{s})$ space. It is to this end that we define the concept of layers. 

Let $\lambda_{u}$ and $\lambda_{s}$ be the multipliers for the internal dynamics on $W_{1}^{u}$ and $W_{2}^{s}$, respectively. Let $\mathbb{T} \times (-D_{u}, D_{u})$ and $\mathbb{T} \times (-D_{s}, D_{s})$ be fundamental domains of the parameterizations of $W_{1}^{u}$ and $W_{2}^{s}$, respectively.  Now, define subsets $U_{n}^{+}$, $U_{n}^{-}$  and $S_{n}^{+}$, $S_{n}^{-}$ of $W_{1}^{u}$ and $W_{2}^{s}$ as follows:
\begin{gather} U_{n}^{+} = \{ W_{1}^{u}(\theta,s) : (\theta,s) \in \mathbb{T} \times [D_{u}\lambda_{u}^{n-1}, D_{u}\lambda_{u}^{n}]\} \\
U_{n}^{-} = \{ W_{1}^{u}(\theta,s) : (\theta,s) \in \mathbb{T} \times [-D_{u}\lambda_{u}^{n-1}, -D_{u}\lambda_{u}^{n}]\} \\
S_{n}^{+} = \{ W_{2}^{s}(\theta,s) : (\theta,s) \in \mathbb{T} \times [D_{s}/\lambda_{s}^{n-1}, D_{s}/\lambda_{s}^{n}]\} \\
S_{n}^{-} = \{ W_{2}^{s}(\theta,s) : (\theta,s) \in \mathbb{T} \times [-D_{s}/\lambda_{s}^{n-1}, -D_{s}/\lambda_{s}^{n}]\} \end{gather}
where $n \in \mathbb{Z}$. Finally, define $U_{n} = U_{n}^{+} \cup U_{n}^{-}$ and $S_{n} = S_{n}^{+} \cup S_{n}^{-}$. We refer to the subsets $U_{n}$ and $S_{n}$ as layers, and to $U_{n}^{+} $, $S_{n}^{+} $ and $U_{n}^{-}$, $S_{n}^{-}$ as positive and negative half-layers, respectively. In our experience, $W_{1}^{u}(\theta_{u}, s_{u})$ and $W_{2}^{s}(\theta_{s}, s_{s})$ do not intersect for $|s_{u}| < D_{u}$ and  $|s_{s}| < D_{s}$; this can usually be seen from plotting the projections of the manifolds for these $s$-values in $(x,y,p_{x})$ space.  Hence, if $W_{1}^{u}$ and $W_{2}^{s}$ intersect, it must be that $U_{n_{1}}$ intersects  $S_{n_{2}}$ for some $n_{1}, n_{2} \in \mathbb{Z}^{+}$. 

The most important property of these layers is that due to Equation \eqref{invariancequationfinal}, $F(U_{n}) = U_{n+1}$ and $F(S_{n}) = S_{n-1}$; more generally, $F^{k}(U_{n}) = U_{n+k}$ and $F^{k}(S_{n}) = S_{n-k}$ for all $k \in \mathbb{Z}$. This allows us to restrict our heteroclinic connection search to only certain pairs of layers of $W_{1}^{u}$ and $W_{2}^{s}$. To see this, suppose we are searching for a heteroclinic connection which comes from layer $U_{n_{1}}$ intersecting layer $S_{n_{2}}$ at $\bold{x} \in \mathbb{R}^{4}$. Then, since $F(U_{n}) = U_{n+1}$ and $F(S_{n}) = S_{n-1}$, we have that $F(\bold{x})$ must be in both $U_{n_{1}+1}$ and $S_{n_{2}-1}$. More generally, for all $k \in \mathbb{Z}$, we have that 
\begin{equation} \label{layerintersection} F^{k}(\bold{x}) \in U_{n_{1}+k} \cap S_{n_{2}-k} \end{equation} 
Now, if $n_{1}$ and $n_{2}$ are both odd or both even, using $k = \frac{n_{2} - n_{1}}{2} $ in Equation \eqref{layerintersection} gives us $F^{k}(\bold{x}) \in U_{\tilde n} \cap S_{\tilde n}$, where $ \tilde n \stackrel{\text{def}}{=} \frac{n_{1} + n_{2}}{2}$. On the other hand, if $n_{1}$ and $n_{2}$ are of opposite parity, setting $k = \frac{n_{2} - n_{1}+1}{2}$ in Equation \eqref{layerintersection} gives us $F^{k}(\bold{x}) \in U_{\tilde n} \cap S_{\tilde n -1}$, where $ \tilde n \stackrel{\text{def}}{=} \frac{n_{1} + n_{2} + 1}{2}$. 

When searching for the heteroclinic trajectory which arises due to the manifolds' intersection at $\bold{x}$, it is enough to find any point on the orbit of $\bold{x}$ under the map $F$, including $F^{k}(\bold{x})$ from the preceding analysis. Based on the above discussion, it is clear that we will find the point $F^{k}(\bold{x})$ if we look for intersections of pairs of layers of form $(U_{n}, S_{n})$ or $(U_{n}, S_{n-1})$ for $n \in \mathbb{Z^{+}}$ (as mentioned earlier, our experience is that the manifolds do not intersect for $|s_{u}| < D_{u}$ and  $|s_{s}| < D_{s}$, so we only consider positive $n$). Since $\bold{x}$ was an arbitrary heteroclinic point, if we search for intersections of pairs of layers of the form just presented above, we will find all possible heteroclinic trajectories. 

As a final note, it is easy to see that if $U_{n_{1}}$ intersects $S_{n_{2}}$, then the time of flight of the resulting heteroclinic connection from the fundamental domain of one torus manifold to the other is $2\pi(n_{1}+n_{2})/\Omega_p$; this is because $n_{1}+n_{2}$ mappings by $F$ are required. Hence, the layer indices can be thought of as a proxy for the connection trajectory time of flight. 

\section{Rapid GPU-Assisted Search for Manifold Intersections} \label{gpuSearchSection}

With methods of computing manifolds and restricting the connection search to certain layer pairs now developed, we next seek to develop computationally fast methods of finding intersections of the manifold layers. The manifolds being dealt with are 2D geometric objects in 4D space, so by discretizing the manifolds and applying methods inspired by those from computer graphics collision detection algorithms, we are able to very rapidly search a pair of manifolds for intersections. The algorithms are massively sped up by taking advantage of the huge number of threads available on modern GPUs. We start this section a with description of the manifold discretization. We then give the full explanation of our algorithm, and finally demonstrate an example application. 

\subsection{Discrete Mesh Representation of Manifolds} \label{meshSection}

In Section \ref{backgroundsection}, we described how it is possible to compute Fourier-Taylor parameterizations $W(\theta, s)$ of the stable/unstable manifolds of stroboscopic map invariant circles. We also gave Equations \eqref{globou} and \eqref{globos} demonstrating how to use the parameterizations to compute $W(\theta, s)$ for $s$-values outside the fundamental domain $\mathbb{T} \times (-D,D)$. Now, we wish to use these tools to generate a discrete representation of the globalized manifold that can be used for computations and analysis. 

When we compute functions of $\theta$, such as $K(\theta)$, $P(\theta)$, $\Lambda(\theta)$, or the manifold $s^{k}$ coefficients $W_k(\theta)$, we represent them on the computer as arrays of function values at $N$ evenly spaced $\theta$ values $\theta_{i} = 2\pi i/N$, $i= 0, 1, \dots, N-1$. By also taking a grid of $2K+1$ evenly spaced $s$-values $s_{k} = kD/K$, $k=-K, -K+1, \dots, 0, 1, \dots, K$ from $-D$ to $D$, we end up with a set of $N(2K+1)$ ordered pairs $(\theta_{i}, s_{k})$. As all these pairs belong to the fundamental domain, we can simply evaluate our Fourier-Taylor parameterizations to compute and store $W(\theta_{i}, s_{k})$ for each $i = 0, \dots, N-1$ and $k = -K,  \dots, K$, getting a set of $N(2K+1)$ points on the manifold. 

After the initial grids of values $W(\theta_{i}, s_{k})$ have been stored, one must numerically globalize the manifolds. We describe the case of an unstable manifold with multiplier $\lambda>1$; the stable case is very similar except for the use of $F^{-1}$ and a few sign changes. To start, fixing $k$ and applying Equation \eqref{globou} to $W(\theta_{i}, s_{k})$ for all $i=0,1,\dots, N-1$ gives $N$ manifold points $F^{n}(W(\theta_{i}, s_{k})) = W(\theta_{i}+n\omega, \lambda^{n}s_{k})$. These are all at the same $s = \lambda^{n}s_{k}$ value but are at $\theta$ values shifted from the $\theta_{i}$. We want all of our manifold points to have the same $\theta$ values; hence, to shift them back, we use a fast Fourier transform (FFT) based translation algorithm. 

Suppose one has a periodic function $a(\theta)$ with known values at $\theta_{i} = 2\pi \frac{i}{N}$, $i= 0, 1, \dots, N-1$. Then, we can take the FFT of this array of function values to find the first $N$ Fourier coefficients $\hat a(i)$, $i= 0, 1, \dots, N-1$. Finally, using the usual formula for $a(\theta_{i})$ given $\hat a(i)$ gives
\begin{equation} \label{translate} a(\theta_{i})  = \frac{1}{N} \sum_{k =0}^{N-1} \hat a(k) e^{jk\theta_{i}}  \rightarrow a(\theta_{i}+\rho)  = \frac{1}{N} \sum_{k =0}^{N-1} [\hat a(k)e^{jk\rho}] e^{jk\theta_{i}}  \end{equation}
where $j$ denotes $\sqrt{-1}$ in Equation \eqref{translate}. Hence, given the $a(\theta_{i})$ values, to find the values $a(\theta_{i}+\rho) $, one takes the FFT, multiplies the $k$th Fourier coefficient by $e^{jk\rho}$, and takes the inverse FFT. By using this algorithm with $\rho=-n\omega$, given the $N$ values $W(\theta_{i}+n\omega, \lambda^{n}s_{k})$ found earlier, we can find $W(\theta_{i}, \lambda^{n}s_{k})$. We now compute and store the points $W(\theta_{i}, \lambda^{n}s_{k})$, along with the corresponding $s$-values, for $n = 1, \dots, n_{max}$, up to some $n_{max} \in \mathbb{Z}^{+}$. We do this for all $k = -K,  \dots, K$. 

The $s$ values $\lambda^{n}s_{k}$ for $n = 0, 1, \dots, n_{max}$, $k = -K, \dots, K$ for the computed points will form an unevenly spaced finite set ranging from $s = -\lambda^{n_{max}}D$ to $\lambda^{n_{max}}D$. Redefine $\{s_{k}\}$ now to be the set of all these $s$ values; for easier notation, we sort and reindex $\{s_{k}\}$ so that $k = 1, \dots, M$ where $M$ is the length of $\{s_{k}\}$. Note that since $-D$ and $D$ were part of our initial grid of $s$ values, the $U_{n}$ layer boundaries (corresponding to points with $s = \pm \lambda^{n}D$) will be contained in our set $\{s_{k}\}$; this fact will be useful later. In our case, we stored the $x, y, p_{x}$, and $p_{y}$ values of the computed manifold points in 4 separate 2D $N \times M$ arrays on the computer, so that the $(i+1,k)$ entry of each array is the $x, y, p_{x}$, or $p_{y}$ coordinate of  $W(\theta_{i}, s_{k})$. Moving down a column of each array corresponds to increasing $\theta$ and constant $s$, and moving across a row corresponds to constant $\theta$ and increasing $s$. 

\begin{figure}
\begin{centering}
\includegraphics[width=0.475\columnwidth]{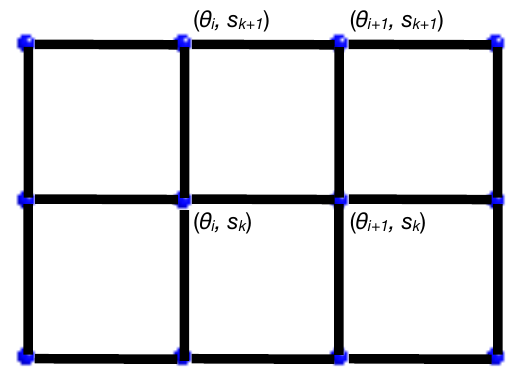}
\includegraphics[width=0.46\columnwidth]{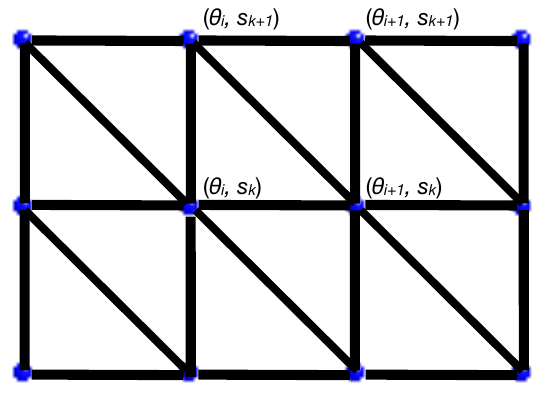}
\caption{ \label{fig:meshschematic} Schematic of Quadrilateral (Quad) and Triangular Mesh Construction}
\end{centering}
\end{figure}

With manifold points computed on a discrete grid of $(\theta, s)$ ordered pairs, we have the numerical results required to obtain the mesh representation of the manifold. We have the values of $W(\theta, s)$ at the points $\{(\theta_{i}, s_{k})\}$ for $i = 0, \dots, N-1$, $k = 1, \dots, M$. Consider the index $i$ to be modulo $N$, so that $\theta_{N} = \theta_{0}$ and $\theta_{-1} = \theta_{N-1}$. To form the manifold mesh, connect $W(\theta_{i}, s_{k})$ with $W(\theta_{i-1}, s_{k})$, $W(\theta_{i+1}, s_{k})$, $W(\theta_{i}, s_{k-1})$, and $W(\theta_{i}, s_{k+1})$ using line segments. If $k-1$ or $k+1$ is outside the range of allowed indices $1, \dots, M$ (which is true if $k=1$ or $k=M$, respectively), then omit the corresponding segment from the mesh. This yields a quadrilateral mesh representation of $W$, as is schematically illustrated on the left of Figure \ref{fig:meshschematic}. We denote the $(i,k)$ quadrilateral (quad, for short) to be that with vertex set $Q_{ik} = \{ W(\theta_{i}, s_{k}), W(\theta_{i+1}, s_{k}), W(\theta_{i}, s_{k+1}), W(\theta_{i+1}, s_{k+1})\} $, where $i=0, \dots, N-1$, $k=1, \dots, M-1$ enumerate the $N(M-1)$ quads in the mesh.

As the vertices of a quad in 4D do not determine a plane, it is better to consider each quad as being composed of two triangles to allow for linear computations with the mesh. We split the $(i,k)$ quad into two triangles by connecting the vertex $W(\theta_{i+1}, s_{k})$ with $W(\theta_{i}, s_{k+1})$. Hence, for each ordered pair $(i,k)$, $i=0, \dots, N-1$, $k=1, \dots, M-1$, we have two triangles. This gives a triangular mesh for $W$. The right of Figure \ref{fig:meshschematic} shows a schematic representation of the mesh construction, illustrating all the points which are connected to each other. 

\begin{figure}
\begin{centering}
\includegraphics[width=0.5\columnwidth]{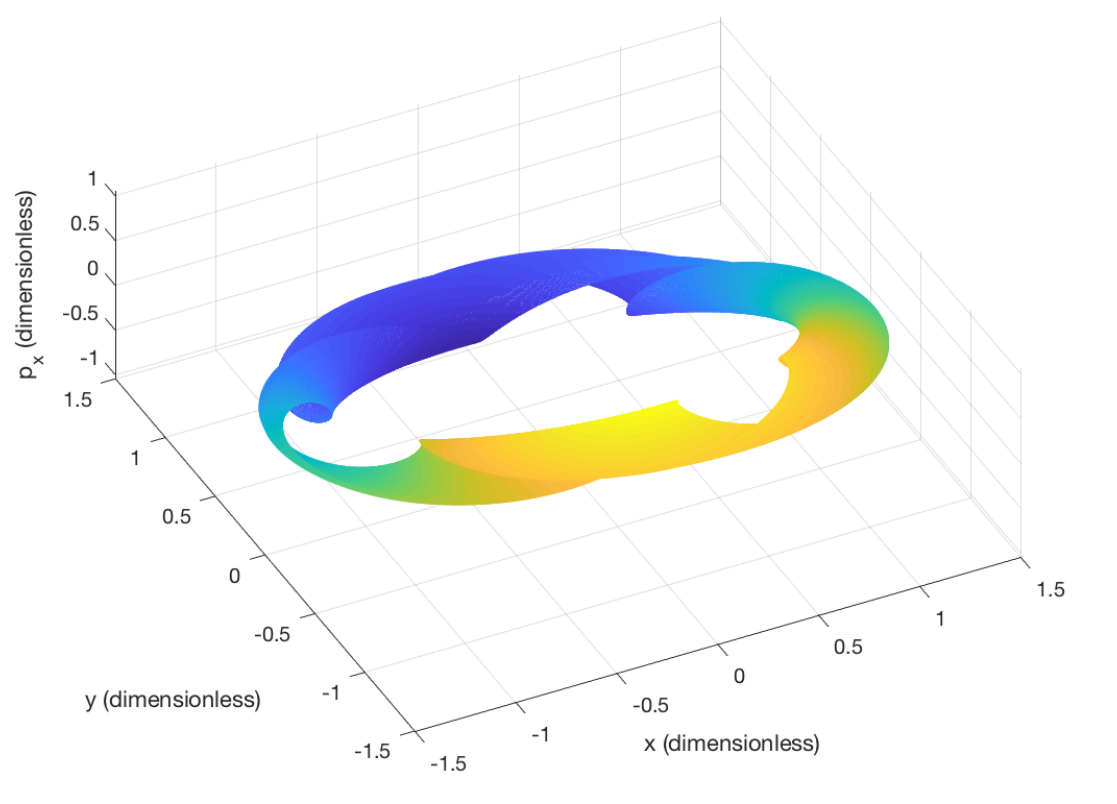}
\caption{ \label{fig:globalmani} $(x,y,p_{x})$ projection of 3:4 $W^{s}$ in Jupiter-Europa PERTBP for $\omega=1.559620297$ (colors added for visualization purposes only)}
\end{centering}
\end{figure}

A 3D projection of an example globalized stable manifold (denoted $W^{s}$) of a 3:4 resonant Jupiter-Europa PERTBP invariant circle is given in Figure \ref{fig:globalmani}; this figure was generated using MATLAB's mesh function, which generates a quad mesh similar to the one described here. 

\subsection{GPU-Accelerated Manifold Mesh Intersection Search}

Now that we have described how to construct the quadrilateral and triangular mesh representations of the manifolds, we start searching for heteroclinic connections. As defined earlier, let $W^{u}_{1}(\theta_{u}, s_{u})$ and $W^{s}_{2}(\theta_{s}, s_{s})$ represent the unstable and stable manifolds of stroboscopic map invariant circles 1 and 2, respectively. Let $\lambda_{u}$ and $\lambda_{s}$ be the multipliers for their internal dynamics and $\mathbb{T} \times (-D_{u}, D_{u})$ and $\mathbb{T} \times (-D_{s}, D_{s})$ be the fundamental domains of the parameterizations for $W_{1}^{u}$ and $W_{2}^{s}$, respectively. After computing and storing the vertices of the meshes for $W_{1}^{u}$ and $W_{2}^{s}$, the problem of finding heteroclinic connections becomes that of finding intersections of these two 2D meshes in 4D space. 

From Section \ref{layerSection}, we know that we can restrict our attention to finding intersections of only certain layers of the manifolds; using the notation from earlier, we seek to find intersections of pairs of layers of the form ($U_{n}$, $S_{n}$) or ($U_{n}$, $S_{n-1}$). Equivalently, we need to check if $U_{n}^{+}$ or $U_{n}^{-}$ intersect any of $S_{n}^{+}$,  $S_{n}^{-}$, $S_{n-1}^{+}$, or $S_{n-1}^{-}$, for $n \in \mathbb{Z}^{+}$. Recalling that the $s$-values generated during manifold globalization include those for the boundaries of these layers and half-layers, it turns out that the half-layers just correspond to easily identified subsets of the manifold meshes. For $U_{n}^{+}$, one simply takes the $W_{1}^{u}$ mesh vertices which satisfy $ s_{u} \in [D_{u} \lambda_{u}^{n-1} , D_{u} \lambda_{u}^{n} ]$. For $S_{n}^{+}$, take vertices of the $W_{2}^{s}$ mesh with $ s_{s} \in [D_{s}/\lambda_{s}^{n-1} , D_{s}/\lambda_{s}^{n} ]$. The negative half-layers are the same except for a change in the signs of $D_{u}$ and $D_{s}$.  If the manifold coordinates are stored in four 2D arrays as described in Section \ref{meshSection}, with each column containing the coordinates of all points for a given $s$ value, then the vertex set of a half-layer mesh is just comprised of points with coordinates from a contiguous set of columns. 

With the meshes for the half-layers identified, we finally arrive at the problem of searching for intersections of two half-layer meshes. Let us say that the mesh vertices corresponding to an unstable half-layer are $W_{1}^{u}(\theta_{u,i}, s_{u,k})$, $i=0, \dots, N_{1}-1$, $k=1, \dots, M_{1}$; for the stable half-layer mesh let the vertices be $W_{2}^{s}(\theta_{s,j}, s_{s,\ell})$, $j=0, \dots, N_{2}-1$, $\ell = 1, \dots, M_{2}$. As was done in the previous section, define the quad vertex sets
\begin{gather*} Q^{u}_{ik} = \{ W_{1}^{u}(\theta_{u,i}, s_{u,k}), W_{1}^{u}(\theta_{u,i+1}, s_{u,k}), W_{1}^{u}(\theta_{u,i}, s_{u,k+1}), W_{1}^{u}(\theta_{u,i+1}, s_{u,k+1})\} \\
Q^{s}_{j\ell} = \{ W_{2}^{s}(\theta_{s,j}, s_{s,\ell}), W_{2}^{s}(\theta_{s,j+1}, s_{s,\ell}), W_{2}^{s}(\theta_{s,j}, s_{s,\ell+1}), W_{2}^{s}(\theta_{s,j+1}, s_{s,\ell+1})\} \end{gather*} 
where $i=0, \dots, N_{1}-1$; $j=0, \dots, N_{2}-1$; $k=1, \dots, M_{1}-1$; and $\ell = 1, \dots, M_{2}-1$ (again consider the indices $i$ and $j$ to be modulo $N_{1}$ and $N_{2}$, respectively). 

\begin{figure}
\begin{centering}
\includegraphics[width=0.5\columnwidth]{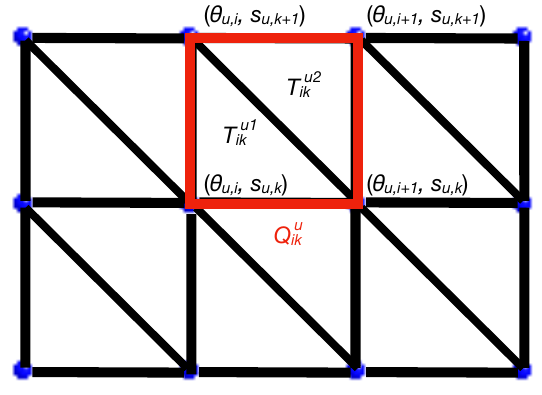}
\caption{ \label{fig:meshnotationschematic} Schematic of Definition of $Q_{ik}^{u}$, $T^{u1}_{ik}$, and $T^{u2}_{ik}$ on Unstable Manifold Mesh}
\end{centering}
\end{figure}

For notational convenience, we also use $Q^{u}_{ik}$ and $Q^{s}_{j\ell}$ to refer to the quads formed by the vertices contained therein. As described earlier, we consider each $Q^{u}_{ik}$ and $Q^{s}_{j\ell}$ to be comprised of two triangles, which we will define by vertex sets
\begin{gather*} T^{u1}_{ik} = \{ W_{1}^{u}(\theta_{u,i}, s_{u,k}), W_{1}^{u}(\theta_{u,i+1}, s_{u,k}), W_{1}^{u}(\theta_{u,i}, s_{u,k+1})\} \\
T^{u2}_{ik} = \{ W_{1}^{u}(\theta_{u,i+1}, s_{u,k}), W_{1}^{u}(\theta_{u,i}, s_{u,k+1}), W_{1}^{u}(\theta_{u,i+1}, s_{u,k+1})\} \\
 T^{s1}_{j\ell} = \{ W_{2}^{s}(\theta_{s,j}, s_{s,\ell}), W_{2}^{s}(\theta_{s,j+1}, s_{s,\ell}), W_{2}^{s}(\theta_{s,j}, s_{s,\ell+1})\} \\
T^{s2}_{j\ell} = \{ W_{2}^{s}(\theta_{s,j+1}, s_{s,\ell}), W_{2}^{s}(\theta_{s,j}, s_{s,\ell+1}), W_{2}^{s}(\theta_{s,j+1}, s_{s,\ell+1})\} \end{gather*} 
Again for ease of notation, we also use $T^{u1}_{ik}$, $T^{u2}_{ik}$, $T^{s1}_{j\ell}$, and $T^{s2}_{j\ell}$ to refer to the plane triangles formed by the vertices contained therein. We have that $Q_{ik}^{u}=T^{u1}_{ik} \cup T^{u2}_{ik}$ and $Q^{s}_{j\ell} = T^{s1}_{j\ell} \cup T^{s2}_{j\ell}$, as is schematically illustrated for the unstable manifold mesh case in Figure \ref{fig:meshnotationschematic}.

Now, the problem of searching for intersections between the half-layers can be solved by checking whether any quad $Q^{u}_{ik}$ intersects any quad $Q_{j\ell}^{s}$; we say that $Q^{u}_{ik}$ intersects  $Q_{j\ell}^{s}$ if any of the triangles $T^{u1}_{ik}$ or $T^{u2}_{ik}$ intersect either of  $T_{j\ell}^{s1}$ or  $T_{j\ell}^{s2}$. This is just the collision detection problem described in Section \ref{colDetectSection} with $A$ as the set of all $Q^{u}_{ik}$ and $B$ being the set of all $Q_{j\ell}^{s}$. There are $N_1 (M_{1}-1)$ quads in the unstable manifold half-layer mesh, and $N_{2} (M_{2}-1)$ in the stable manifold half-layer mesh; hence, $N_1 N_{2} (M_{1}-1) (M_{2}-1)$ pairs of quads must be checked for intersection. For an example computation in the PERTBP computing intersections between manifolds of 3:4 and 5:6 resonant invariant circles, we had $N_{1} = 1024$, $N_{2} = 2048$, and $M_{1} = M_{2} = 35$, for a total of 2,424,307,712 pairs of quads; each quad pair intersection test involves 4 triangle-triangle checks, which gives 9,697,230,848 pairs of triangles to check for our example. Given the massive number of identical checks to be done, it is clear that a GPU will be well suited to this application. 

It is possible to exactly determine whether two 2D triangles intersect in 4D by solving a $4\times4$ system of linear equations and checking whether the solution satisfies certain conditions, as will be described in more detail later (see Section \ref{narrowPhase}); 4 such tests are required in each quad-quad intersection test. However, solving a different $4 \times 4$ system for each of billions of triangle pairs would be extremely computationally expensive, and would not be an efficient algorithm to implement even on a GPU. Fortunately, as is the expectation in practical collision detection problems, the vast majority of pairs of quads will not intersect when checking two manifold half-layer meshes for intersection. Hence, as described in Section \ref{colDetectSection}, it is necessary to first have a computationally cheap broad phase algorithm which can quickly reject most of the non-intersecting quad pairs. For this, we first implement a simple method of grouping quads such that only quads from the same group can possibly intersect; this is then followed by two common broad phase non-intersection tests, both of which can be run on the GPU. 

One thing to note is that the quad/triangle meshes are only approximations of the stable and unstable manifolds, with points other than the mesh vertices represented by linear interpolation. Thus, even an exact intersection between two mesh triangles will only yield an approximate intersection of their manifolds. Nevertheless, the meshes allow for computationally rapid identification of such near-intersections, which will provide good initial guesses for a more accurate calculation. This will be described in more detail in Section \ref{refinementSection}. 

\subsubsection{Broad Phase: Uniform Grid Spatial Partitioning} \label{partitionSection}

The first step of the broad phase algorithm is aimed at excluding pairs of quads which are in completely different regions of phase space. The main idea \cite{collisionDetection} is to partition a finite ``world'' which contains all our objects into a uniform grid of boxes. Then, it is clear that only objects (in our case, quads) overlapping a common box can possibly intersect. Thus, for each box one can make two lists, one of the $Q^{u}_{ik}$ and another list of the $Q_{j\ell}^{s}$ which overlap that box; by taking all the pairs of quads having one quad from each list, one gets the set of quad pairs which potentially intersect in that box. The union of these sets over all boxes forms the set of all potentially intersecting quad pairs. 

\begin{figure}
\centering
\includegraphics[width=0.5\columnwidth]{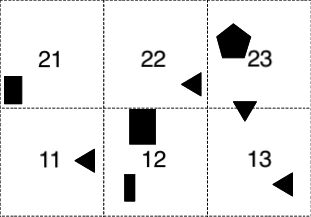}
\caption{ \label{fig:gridpart} Illustration\cite{gridFig} of uniform grid spatial partitioning in 2D space}
\end{figure}

As long as the grid was not too coarse, one will end up with a significantly smaller list of quad pairs after this procedure. This is schematically illustrated in 2D in Figure \ref{fig:gridpart} with polygons instead of quads, and $A=B$ (in the framework of Section \ref{colDetectSection}); here, there are 8 objects, so $28$ pairs of objects in the world. The world is partitioned into a uniform grid of 6 boxes, each identified by two grid indices ranging from 1 to 3 in $x$ and 1 to 2 in $y$. After this, there are only 3 potentially intersecting pairs of objects, one from each of the boxes 12, 13, and 23. Other than cutting the number of quad pairs to check, the key advantage to the spatial partitioning procedure is that finding the lists of quads in each box has complexity $O(|A| + |B|)$ rather than $O(|A||B|)$. To see this, we briefly describe the steps involved. 

First of all, we define a finite ``world'' as a large box containing the quads, and a uniform grid. The world's minimum and maximum $x$-bounds $x_{min}$ and $x_{max}$ can simply be taken as the minimum and maximum, respectively, of $x$ over all the two manifolds' half-layer points $W_{1}^{u}(\theta_{u,i}, s_{u,k})$ and $W_{2}^{s}(\theta_{s,j}, s_{s,\ell})$ (more efficient bounds can be found as well, see remark at the end of this section). The grid size in $x$ can then be set as $\Delta_{x} = \frac{x_{max}-x_{min}}{N_{x}}$ for some $N_{x} \in \mathbb{Z}^{+}$; one should choose $N_{x}$ so that $\Delta_{x}$ is greater than the largest $x$-width of all quads. $y, p_{x}$, and $p_{y}$ world and grid sizes are done similarly. Then, for each quad, one calculates its overlapped grid indices in each coordinate; for example, given $Q^{u}_{ik}$ with minimum and maximum $x$-coordinates $x_{min, ik}$ and $x_{max, ik}$, the overlapped $x$ grid indices are $\lceil \frac{x_{min,ik}-x_{min}}{\Delta_{x}} \rceil$ and $\lceil \frac{x_{max,ik}-x_{min}}{\Delta_{x}} \rceil$. The same can be done with the $Q^{s}_{j\ell}$ as well as with the $y, p_{x}$, and $p_{y}$ grid indices. 

Once the overlapped grid indices in $x, y, p_{x}$, and $p_{y}$ have been found for all quads, one forms the lists of $Q^{u}_{ik}$ and $Q_{j\ell}^{s}$ overlapping each box. For this, one can simply iterate over all boxes, each of which corresponds to a tuple of grid indices, and find all $Q^{u}_{ik}$ and $Q^{s}_{j\ell}$ overlapping that same combination of indices. Functions such as MATLAB's {\tt find} or Julia's {\tt findall}, which return all true indices of an array, can be useful for this step. It is clear that at no point of this algorithm are pairs of quads considered; all computations involve only one quad at a time, hence the $O(|A| + |B|)$ complexity rather than $O(|A||B|)$. Finally, from these lists one finds all potentially intersecting quad pairs in each box and in the world as a whole. 

\begin{remark}

One can also define the world bounds in such a manner that excludes quads which have no possibility of intersecting any quad on the other manifold. We show how to do this for the lower $x$-bound $x_{min}$; the upper $x$ as well as all $y$, $p_{x}$, and $p_{y}$ bounds are found similarly. First, one should find the minimum $x$-value for all the half-layer points $W_{1}^{u}(\theta_{u,i}, s_{u,k})$, as well as the minimum $x$-value for all the half-layer points $W_{2}^{s}(\theta_{s,j}, s_{s,\ell})$, and then let $x_{min,max}$ be the larger of these two values. Now, if $x_{min,max}$ was the minimum of the $W_{1}^{u}(\theta_{u,i}, s_{u,k})$, set $x_{min}$ to be the minimum $x$-value across all the quads $Q_{j\ell}^{s}$ which have a maximum $x$-value greater than $x_{min,max}$ (this step is necessary in order not to exclude quads on $W_{2}^{s}$ which might intersect the minimum $x$-value quad of $W_{1}^{u}$). Or, if $x_{min,max}$ was the minimum of the $W_{2}^{s}(\theta_{s,j}, s_{s,\ell})$, simply reverse the roles of manifolds 1 and 2 in the preceding step. This gives a lower bound $x_{min}$ such that when the grid indices are computed, one can ignore any quad whose overlapped grid indices are all less than 1. In fact, when trying to set $x_{min}$ using this algorithm, one may find in the second step that no quad on the other manifold has a maximum $x$-value greater than $x_{min,max}$, in which case one immediately concludes that the half-layers cannot intersect. 

\end{remark}

\subsubsection{Broad Phase: Bounding Box Test} \label{bdboxSection}

\begin{figure}
\includegraphics[width=0.495\columnwidth]{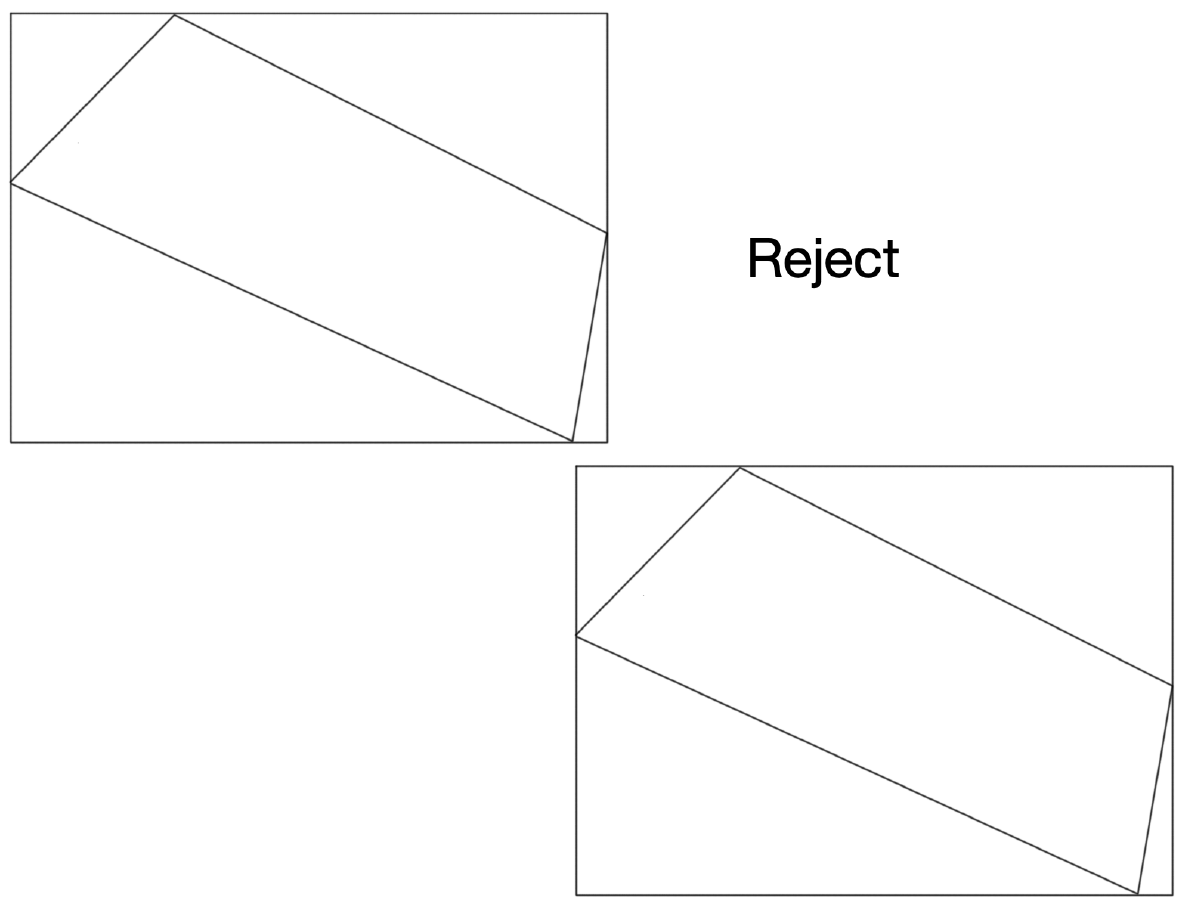}
\includegraphics[width=0.495\columnwidth]{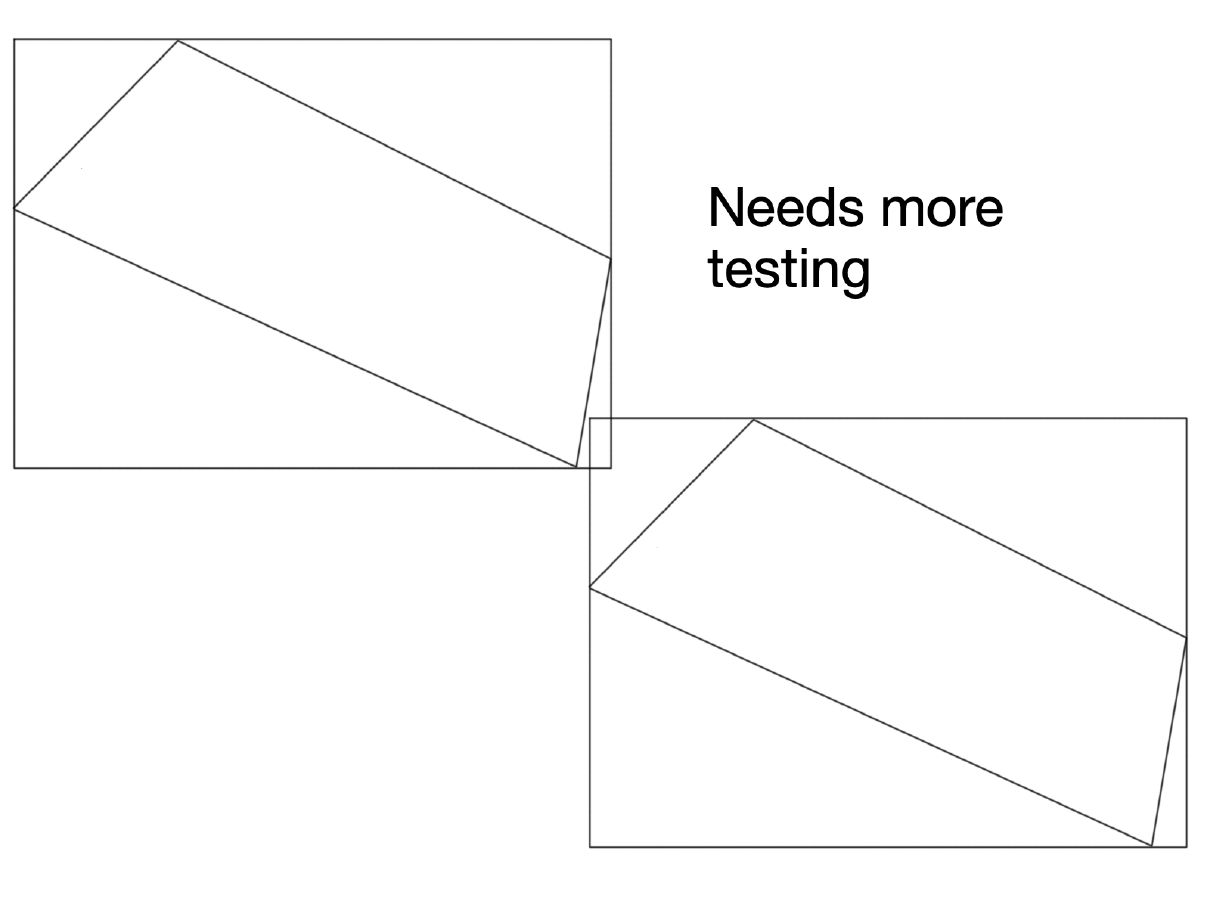}
\caption{ \label{fig:bdbox} Illustration of bounding box test in 2D space}
\end{figure}

The next step of the broad phase is a simple axis-aligned bounding box test\cite{collisionDetection}, applied to each of the potentially intersecting quad pairs found in the previous step. This test is generally used with geometric objects in 3D space, but the concept works the same in 4D. The basic idea is to consider each quad to be enclosed by a minimal 4D box (the bounding box) having its edges parallel to the $x, y, p_{x}$, and $p_{y}$ axes. Then, to check whether two quads $Q_{ik}^{u}$ and $Q_{j\ell}^{s}$ might intersect, simply check whether their corresponding bounding boxes intersect. If they do not, then we can reject the possibility of the quads intersecting; if the boxes do intersect, then additional testing is required. Figure \ref{fig:bdbox} illustrates how this test works in 2D. The test is equivalent to checking whether the maximum $x$ coordinate of the 4 vertices of $Q_{ik}^{u}$ is less than the minimum $x$ coordinate of the 4 vertices of $Q_{j\ell}^{s}$; we also reverse the roles of $Q_{ik}^{u}$ and $Q_{j\ell}^{s}$ and repeat the check. The same is also done for the $y$, $p_{x}$, and $p_{y}$ coordinates; if any of the checks are true, then the quads cannot intersect. 

\subsubsection{Broad Phase: M\"oller Quick Test} \label{mollerSection}

Even after the spatial grid partitioning and bounding box tests have excluded most of the non-intersecting quad pairs, there are still often a fair number of pairs left. Hence, we run a second broad phase test on the remaining pairs which is due to M\"oller \cite{moller}. The full M\"oller test is a two-part triangle-triangle intersection test developed for applications in 3D graphics; we use only the first part of this test, which is for quick rejection of non-intersecting pairs of triangles. Given two triangles in 3D space, the quick test checks whether all the vertices of one triangle are on the same side of the plane formed by the other; if so, the pair of triangles cannot intersect. In our case, since our objects are in 4D $(x,y,p_{x},p_{y})$ space, we project the quads and their constituent triangles onto 3D $(x,y,p_{x})$ space in order to carry out the test. Of course, if the 3D projected quads and triangles do not intersect, then neither will the full quads and triangles in 4D space.  

Denote the projected quads in $(x,y,p_{x})$ space as $\tilde Q_{ik}^{u}$ and $\tilde Q_{j\ell}^{s}$, and their constituent projected triangles as $\tilde T^{u1}_{ik}$, $\tilde T^{u2}_{ik}$ and $\tilde T_{j\ell}^{s1}$, $\tilde T_{j\ell}^{s2}$. Then, the first step in the M\"oller quick test is to see whether all four vertices of $\tilde Q_{j\ell}^{s}$ are (1) on the same side of $\tilde T^{u1}_{ik}$, and (2) on the same side of $\tilde T^{u2}_{ik}$. Both statements (1) and (2) must be true in order to rule out an intersection with $\tilde Q^{u}_{ik}$. To check (1), we first need the equation of the plane in which $\tilde T^{u1}_{ik}$ lies. This can be found using standard techniques; the plane will be comprised of all points $\bold{X} \in \mathbb{R}^{3}$ such that
\begin{equation} \label{planeequation} f_{ik}^{u1}(\bold{X}) \stackrel{\text{def}}{=}  N_{ik}^{u1} \cdot (\bold{X} - \tilde W_{1}^{u}(\theta_{u,i}, s_{u,k})) = 0\end{equation} 
\begin{equation*} N_{ik}^{u1} = \left[\tilde W_{1}^{u}(\theta_{u,i+1}, s_{u,k})-\tilde W_{1}^{u}(\theta_{u,i}, s_{u,k})\right] \times \left[ \tilde W_{1}^{u}(\theta_{u,i}, s_{u,k+1}) - \tilde W_{1}^{u}(\theta_{u,i}, s_{u,k}) \right] \end{equation*} 
where $ \tilde W_{1}^{u}(\theta_u, s_u)$ denotes the projection of the point $ W_{1}^{u}(\theta_u, s_u)$ into $(x,y,p_{x})$ space. With the Equation \eqref{planeequation} of the plane of $\tilde T_{ik}^{u1}$ found, we can determine whether all vertices of $\tilde Q_{j\ell}^{s}$ are on the same side of $\tilde T^{u1}_{ik}$ by simply evaluating $f_{ik}^{u1}(\bold{X})$ at the four vertices and seeing if the resulting four values all have the same sign. If they do, then (1) is true, otherwise it is not. 

To check (2), we do the same procedure, but using the vertices of $\tilde T^{u2}_{ik}$ instead of those of $\tilde T^{u1}_{ik}$. If both (1) and (2) are true, then we reject the possibility of $Q_{j\ell}^{s}$ intersecting $ Q^{u}_{ik}$. Otherwise, we carry out the same test as above, but with the roles of $Q_{j\ell}^{s}$ and $\tilde Q^{u}_{ik}$ swapped so that we check whether all four vertices of $\tilde Q_{ik}^{u}$ are (1) on the same side of $\tilde T^{s1}_{j\ell}$, and (2) on the same side of $\tilde T^{s2}_{j\ell}$. If both of these statements are true, then we conclude that $Q_{j\ell}^{s}$ cannot intersect $ Q^{u}_{ik}$; otherwise, we move on to the final, precise test. 

\subsubsection{Narrow Phase} \label{narrowPhase}

As described in Section \ref{colDetectSection}, the broad phase of collision detection is followed by the narrow phase. Thankfully, after the broad phase algorithm is run on the quad pairs generated from two half-layer meshes, the number of potentially intersecting quad pairs left to test is usually quite small. For our earlier example with 2,424,307,712 pairs of quads to be checked for each pair of half-layers, after the broad phase, there would only be at very most a few tens of thousands of cases left to test, and usually far less. In fact, for pairs of half-layers with smaller indices (such as $U_{2}^{+}$, $S_{2}^{+}$), our experience is that the bounding box test alone rejects all of the quad pairs! For those few pairs of quads which have not been rejected in the broad phase, we now need to run a more computationally heavy narrow phase test to check for possible intersections, as well as computing the intersection if it exists. 

Suppose that $ Q^{u}_{ik}$ and $Q_{j\ell}^{s}$ are such a pair of quads which passed the broad phase. At this stage, we start dealing exclusively with their constituent triangles; we check whether any of $T^{u1}_{ik}$ or $T^{u2}_{ik}$ intersect either of  $T_{j\ell}^{s1}$ or  $T_{j\ell}^{s2}$. For this, we need an algorithm to determine whether (and if so, where) two triangles intersect each other in 4D space.  Let $T_{1}$ and $T_{2}$ be two triangles with vertices $\bold{x_{1}}, \bold{x_{2}}, \bold{x_{3}} \in \mathbb{R}^{4}$ and $\bold{y_{1}}, \bold{y_{2}}, \bold{y_{3}} \in \mathbb{R}^{4}$, respectively. Then, to determine whether $T_{1}$ and $T_{2}$ intersect, the first step is to find the intersection of the planes containing these two triangles. The equation to solve to help find this is 
\begin{equation} \label{precise} \bold{x_{2}} + (\bold{x_{1}}-\bold{x_{2}})a + (\bold{x_{3}}-\bold{x_{2}})b = \bold{y_{2}} + (\bold{y_{1}}-\bold{y_{2}})c + (\bold{y_{3}}-\bold{y_{2}})d \end{equation}
where $a,b,c,d \in \mathbb{R}$ are the quantities to be solved for. Equation \eqref{precise} is a 4D linear equation with 4 unknowns, so it generically admits a unique solution. After solving for $a,b,c,$ and $d$, it is easy to see that the LHS and RHS of Equation \eqref{precise} are themselves equal to the intersection point of the planes containing the two triangles; the values of $a,b,c,$ and $d$ determine whether this intersection point lies inside each triangle itself. The conditions for the intersection to be in $T_{1}$ and $T_{2}$ are simple; we just need $a,b,c,d \geq 0$, $a+b \leq 1$, and $c+d \leq 1$. 

If the three conditions just given are not all satisfied, then $T_{1}$ and $T_{2}$ do not intersect. If they are all satisfied, then it is confirmed that the triangles do intersect and the intersection point is given by either side of Equation \eqref{precise}. Furthermore, if $T_{1}$ and $T_{2}$ are $T^{u1}_{ik}$ or $T^{u2}_{ik}$ and $T_{j\ell}^{s1}$ or  $T_{j\ell}^{s2}$, respectively, we can actually use the values of $a,b,c,$ and $d$ to estimate a solution for the heteroclinic connection equation given at the beginning of Section \ref{layerSection}, repeated below:
\begin{equation} \label{toestimate} W^{u}_{1}(\theta_{u}, s_u) = W^{s}_{2}(\theta_{s}, s_{s}) \end{equation}
If $T_{1} = T^{u1}_{ik}$  with $\bold{x}_{1} = W_{1}^{u}(\theta_{u,i+1}, s_{u,k})$, $\bold{x}_{2} = W_{1}^{u}(\theta_{u,i}, s_{u,k})$, $\bold{x}_{3} = W_{1}^{u}(\theta_{u,i}, s_{u,k+1})$, we have
\begin{equation} (\theta_{u}, s_{u}) \approx ((1-a)\theta_{u,i}+a\theta_{u,i+1}, (1-b)s_{u,k}+bs_{u,k+1} ) \end{equation}
If instead $T_{1} = T^{u2}_{ik}$  with $\bold{x}_{1} = W_{1}^{u}(\theta_{u,i}, s_{u,k+1})$, $\bold{x}_{2} = W_{1}^{u}(\theta_{u,i+1}, s_{u,k+1})$, $\bold{x}_{3} = W_{1}^{u}(\theta_{u,i+1}, s_{u,k})$, 
\begin{equation} (\theta_{u}, s_{u}) \approx (a\theta_{u,i}+(1-a)\theta_{u,i+1}, bs_{u,k}+(1-b)s_{u,k+1} ) \end{equation}
Similarly, if $T_{2} = T^{s1}_{j\ell}$, $\bold{y}_{1} = W_{2}^{s}(\theta_{s,j+1}, s_{s,\ell})$, $\bold{y}_{2} = W_{2}^{s}(\theta_{s,j}, s_{s,\ell})$, $\bold{y}_{3} = W_{2}^{s}(\theta_{s,j}, s_{s,\ell+1})$; or $T_{2} = T^{s2}_{j\ell}$, $\bold{y}_{1} = W_{2}^{s}(\theta_{s,j}, s_{s,\ell+1})$, $\bold{y}_{2} = W_{2}^{s}(\theta_{s,j+1}, s_{s,\ell+1})$, $\bold{y}_{3} = W_{2}^{s}(\theta_{s,j+1}, s_{s,\ell})$, then respectively 
\begin{equation} (\theta_{s}, s_{s}) \approx ((1-c)\theta_{s,j}+c\theta_{s,j+1}, (1-d)s_{s,\ell}+ds_{s,\ell+1}) \end{equation}
\begin{equation} (\theta_{s}, s_{s}) \approx (c\theta_{s,j}+(1-c)\theta_{s,j+1}, ds_{s,\ell}+(1-d)s_{s,\ell+1}) \end{equation}
We then store these approximate solutions $(\theta_{u}, s_u, \theta_{s}, s_{s}) $ of Equation \eqref{toestimate}, as well as the corresponding manifold mesh intersection points found by solving Equation \eqref{precise}. 

One thing to note about this search for manifold intersections is that we represent the manifolds using planar meshes, since the mesh faces are 2D triangles. The true manifolds lie close to these triangular faces, but the manifolds are curved rather than planar. Hence, any intersection found from this search will be subject to some error, on the order of the squares of $a,b,c$, and $d$. We will describe how to refine these manifold intersections to higher accuracy later in Section \ref{refinementSection}. Conversely, if the mesh intersection test excludes a pair of triangles from intersecting by a sufficiently large distance, then more accurate nearby manifold representations  will also not intersect. 

\subsection{Computational Implementation} \label{compImplementSection}

With the various steps of our mesh intersection algorithm explained, we now describe the implementation of our manifold mesh intersection method in a computer program. Our programs were written using the Julia programming language\cite{bezanson2017julia}, a relatively new high-level language which has gained significant interest in recent years due to its excellent performance, ease of use, multiple-dispatch features, and abundance of high-quality packages (many of which work seamlessly with each other, thanks to multiple dispatch). There is a Julia package called OpenCL.jl\cite{opencljl} which allows one to transfer data to and from an OpenCL device and run OpenCL kernels from Julia programs; the C code for the kernel is simply passed as a large string to an OpenCL.jl function, which generates a kernel which can be run from within Julia. Identification of manifold mesh vertices belonging to a certain half-layer is just a matter of careful indexing; most of the computations occur in trying to detect intersections of two half-layers, so it is this part of the method we focus on here. 

As was defined earlier, let $W_{1}^{u}(\theta_{u,i}, s_{u,k})$, $i=0, \dots, N_{1}-1$, $k=1, \dots, M_{1}$ and $W_{2}^{s}(\theta_{s,j}, s_{s,\ell})$, $j=0, \dots, N_{2}-1$, $\ell = 1, \dots, M_{2}$ be the vertices of the unstable and stable manifold half-layer meshes being considered, respectively. We store the vertex coordinates of each manifold half-layer mesh in four 2D arrays, one for each of $x,y, p_{x}$, and $p_{y}$;  using the convention of 1-based indexing as in MATLAB and Julia, the $(i+1,k)$ entries of the four $N_{1} \times M_{1}$ arrays containing unstable half-layer coordinates are the $x,y, p_{x}$, and $p_{y}$ coordinates of $W_{1}^{u}(\theta_{u,i}, s_{u,k})$. Similarly, the $(j+1,\ell)$ entries of the four $N_{2} \times M_{2}$ arrays containing stable half-layer coordinates are the $x,y, p_{x}$, and $p_{y}$ coordinates of $W_{2}^{s}(\theta_{u,j}, s_{u,\ell})$. These eight 2D coordinate arrays are provided as inputs to our mesh intersection function, which starts by applying the spatial partitioning scheme of Section \ref{partitionSection} to determine the world, grid, and lists of $Q^{u}_{ik}$ and $Q_{j\ell}^{s}$ overlapping each grid box; we store the lists of quads in two arrays of arrays, one for the $Q^{u}_{ik}$ and another for the $Q_{j\ell}^{s}$. The $m$\textsuperscript{th} element of each array of arrays is a 1D array of 0-based linear indices ($i+(k-1)N_{1}$ for $Q^{u}_{ik}$ or $j+(\ell-1)N_{2}$ for $Q_{j\ell}^{s}$) identifying the quads overlapping the $m$\textsuperscript{th} grid box. Note that a box cannot contain an intersecting pair of quads if no $Q^{u}_{ik}$ or no $Q_{j\ell}^{s}$ overlaps it; hence, we discard all elements of the arrays of arrays corresponding to such boxes. 

The spatial partitioning, since it does not involve pairs and is thus linear in the number of quads, is done on the CPU. However, the remaining broad phase tests apply to the potentially intersecting quad pairs identified by the partitioning; the number of such pairs can still be a fairly large number. Thus the bounding box and quick M\"oller tests of Sections \ref{bdboxSection}-\ref{mollerSection} are done on the GPU using OpenCL.jl. However, OpenCL C versions prior to 2.0, including the 2023 version (1.2) for MacOS, do not support arrays of arrays (pointers to pointers). Hence, we first convert the two arrays of 1D arrays of $Q^{u}_{ik}$, $Q_{j\ell}^{s}$ linear indices into two large 1D arrays of quad indices \texttt{idxs\_quads\_1,2} by simply concatenating the quad index arrays over all grid boxes. The \texttt{idxs\_quads\_1,2} arrays are supplemented by auxiliary arrays \texttt{grid\_box\_nums\_1,2}, whose $m$\textsuperscript{th} entries are the number of $Q^{u}_{ik}$, $Q_{j\ell}^{s}$ overlapping the $m$\textsuperscript{th} box. We also compute the cumulative sum arrays of \texttt{grid\_box\_nums\_1,2}, to each of which we then prepend a 0. We denote the resulting arrays by \texttt{box\_start\_idxs\_1,2}, as the $m$\textsuperscript{th} entry of these two arrays will be the 0-based starting index in \texttt{idxs\_quads\_1,2} of the $m$\textsuperscript{th} box's $Q^{u}_{ik}$ or $Q_{j\ell}^{s}$ indices list. 

The array \texttt{box\_start\_idxs\_1,2} will allow the OpenCL kernel to identify the elements of \texttt{idxs\_quads\_1,2} corresponding to a given grid box. However, we still need a way of determining which box a given thread should investigate. For this we first compute the element-wise product of \texttt{grid\_box\_nums\_1} and \texttt{grid\_box\_nums\_2}, which gives an array whose $m$\textsuperscript{th} entry is the number of potentially intersecting quad pairs in the $m$\textsuperscript{th} box. Again taking its cumulative sum and prepending 0, we store the end result as an array \texttt{box\_gid\_idxs}; we will be able to use this to assign a box to each OpenCL work item (see lines \ref{whileLoop}-\ref{endWhileLoop} of Algorithm \ref{kernelAlg}). Its $m$\textsuperscript{th} entry is the number of potentially intersecting quad pairs in boxes 0 through $m-1$, so its last entry gives the total number $N_{total}$ of such pairs to be considered across all boxes. Finally, we transfer the manifold half-layer coordinate arrays as well as \texttt{idxs\_quads\_1,2},  \texttt{box\_start\_idxs\_1,2}, and \texttt{box\_gid\_idxs} to the GPU. We also allocate an output buffer, denoted \texttt{out}, of $N_{total}$ 32-bit unsigned integers on the GPU, where we will store identifiers for quad pairs which pass all broad phase tests. We also supply the GPU with a pointer \texttt{num\_passed} to an integer initialized to zero, which will serve as a counter for how many quad pairs are not rejected. 

Once the data is transferred to the GPU and buffers allocated, we execute our bounding box and quick M\"oller kernel (recall that the narrow phase in this case is unsuitable for GPU implementation). We need a way to determine which two quads $Q^{u}_{ik}$ and $Q_{j\ell}^{s}$ each work item should test. Hence, the first few steps of the kernel involve extracting the indices of these quads from \texttt{idxs\_quads\_1,2}. As OpenCL kernels are written in a version of C, the 2D indexing used in Julia does not apply here when reading from the arrays of half-layer coordinates; this is the reason that we stored linear indices in the \texttt{idxs\_quads\_1,2} arrays earlier. Each work item applies the broad phase tests to a different pair of quads. The entire kernel is too long to include in this paper, but the essential steps are summarized in Algorithm \ref{kernelAlg}. 
\begin{algorithm}
\caption{Broad phase OpenCL kernel}
\label{kernelAlg}
\begin{algorithmic}[1]
\STATE \texttt{int gid} $\leftarrow$ work item global ID (will range from 0 to $N_{total}-1$) 
\STATE  \texttt{int bid} $\leftarrow$ 0
\WHILE{\texttt{box\_gid\_idxs}[\texttt{bid}+1] $\leq$ \texttt{gid}} \label{whileLoop}
     \STATE \texttt{bid}++ (loop determines which box the work item \texttt{gid} should work with, while making sure enough work items are assigned to each box to process all quad pairs in that box)
\ENDWHILE  \label{endWhileLoop}
\STATE Find starting index in \texttt{idxs\_quads\_1,2} for box \texttt{bid}'s overlapping $Q^{u}_{ik}$, $Q_{j\ell}^{s}$ indices list: \texttt{box\_start\_idx\_1,2}  $\leftarrow$ \texttt{box\_start\_idxs\_1,2}[\texttt{bid}]
 \STATE Find number of quads $Q^{u}_{ik}$, $Q^{s}_{j\ell}$ overlapping box \texttt{bid}: 
 
 \texttt{box\_num\_quads\_1,2}  $\leftarrow$ \texttt{box\_start\_idxs\_1,2}[\texttt{bid}+1]-\texttt{box\_start\_idx\_1,2}

 \STATE Find a linear identifier for which pair of quads in box \texttt{bid} this work item will test:
 
 \texttt{pid\_box}  $\leftarrow$ \texttt{gid} - \texttt{box\_gid\_idxs}[\texttt{\texttt{bid}}] 
      \STATE Convert the linear identifier \texttt{pid\_box} to a pair of quad list indices for box \texttt{bid}. Read the quad lists at those indices to find the corresponding quads' coordinate array linear indices:
      
      \texttt{idx\_1}  $\leftarrow$ \texttt{idxs\_quads\_1}[\texttt{box\_start\_idx\_1} + (\texttt{pid\_box} \% \texttt{box\_num\_quads\_1})]

  \texttt{idx\_2}  $\leftarrow$ \texttt{idxs\_quads\_2}[\texttt{box\_start\_idx\_2} + (\texttt{pid\_box} / \texttt{box\_num\_quads\_1})]
 
\STATE Read the coordinates of the $Q^{u}_{ik}$ and $Q_{j\ell}^{s}$ vertices from the coordinate arrays using \texttt{idx\_1} and \texttt{idx\_2}. Three vertices of $Q^{u}_{ik}$ will be stored at \texttt{idx\_1}, \texttt{idx\_1}+1, and \texttt{idx\_1}+$N_{1}$. The fourth will be at \texttt{idx\_1}+$N_{1}$+1 if \texttt{idx\_1}$\% N_{1} \neq -1$, else at \texttt{idx\_1}$-N_{1}$+1.  Similar for $ Q^{s}_{j\ell}$. 

\STATE Use fmin/fmax to find smallest/largest $x,y,p_{x},p_{y}$ coordinates of vertices of $ Q^{u}_{ik}$ and $Q_{j\ell}^{s}$. Carry out bounding box test; store 0 in \texttt{val} if intersection is rejected, otherwise 1. 

\IF{ \texttt{val} == 1} 
\STATE \label{step6} Compute normals $N^{u1}_{ik}$ and $N^{u2}_{ik}$. Carry out quick M\"oller test to check if $ \tilde Q^{s}_{j\ell}$ vertices are all on the same side of $\tilde T^{u1}_{ik}$ and $\tilde T^{u2}_{ik}$ (OpenCL has functions for cross and dot products). If intersection is rejected, \texttt{val} $\leftarrow 0$. \ENDIF

\IF{ \texttt{val} == 1} 
\STATE Repeat step \ref{step6} with roles of $\tilde  Q^{u}_{ik}$ and $ \tilde Q^{s}_{j\ell}$ reversed. If intersection is rejected, \texttt{val} $\leftarrow 0$. \ENDIF

\IF{ \texttt{val} == 1} 
\STATE \label{writeOut} Each thread with \texttt{val} true gets a unique consecutive value of \texttt{out\_idx}. Store an identifier for the potentially intersecting quad pair $Q^{u}_{ik}$, $Q^{s}_{j\ell}$ in the \texttt{out\_idx} entry of \texttt{out}. 

\texttt{out\_idx} $\leftarrow$ \texttt{atomic\_inc}(\texttt{num\_passed})

\texttt{out}[\texttt{out\_idx}] $\leftarrow$ \texttt{idx\_1} + $N_{1} $*$ M_{1}$*\texttt{idx\_2}
 \ENDIF
\end{algorithmic} \end{algorithm}

Once the broad phase kernel has finished, we read the \texttt{num\_passed} counter value $N_{pass}$ to find out how many potentially intersecting quad pair identifiers there are in the \texttt{out} GPU buffer. We then transfer the first $N_{pass}$ entries from \texttt{out} to the host computer memory, after which we convert each pair identifier \texttt{pid} back to a pair of linear quad indices by inverting the expression for \texttt{out}[\texttt{out\_idx}] from line \ref{writeOut} of Algorithm \ref{kernelAlg}; this yields $\texttt{idx\_1} = \texttt{pid} \% (N_{1}*M_{1})$ (for $Q^{u}_{ik}$)  and $\texttt{idx\_2} = \texttt{pid} / (N_{1}*M_{1})$ (for $Q^{s}_{j\ell}$). Finally, using these indices to find the relevant quads' vertex coordinates, we apply the narrow phase test to each potentially intersecting pair of quads listed in \texttt{out}. This test is carried out on the CPU, and the resulting intersections are saved as well as the corresponding approximate solutions $(\theta_{u}, s_{u}, \theta_{s}, s_{s})$ to Equation \eqref{toestimate}. 

\subsection{Computational Results} \label{benchmarkSection}

The Julia program implementing the previous algorithms was tested on three different consumer-grade machines. Device 1 was a 2017-era laptop with a 2.9GHz quad core Intel i7-7820HQ CPU and an AMD Radeon Pro 560 GPU with 4GB VRAM. Device 2 was a 2019-era laptop with a 2.6GHz six core Intel i7-9750H CPU and an AMD Radeon Pro 5300M GPU with 4GB VRAM. Device 3 was a desktop tower with a 2011-era 3.33GHz six core Intel Xeon W3680 CPU and a 2016-era Radeon RX 480 GPU with 8GB VRAM. The application used for benchmarking the algorithm was the computation of connections between 3:4 resonance $W^{u}$ and 5:6 resonance $W^{s}$ manifolds in the Jupiter-Europa PERTBP, globalized until layers $U_{14}$ and $S_{14}$. We had $N_{1} = 1024$, $N_{2} = 2048$, and $M_{1} = M_{2} = 35$. The 3:4 resonant torus whose $W^{u}$ was computed was the one with $\omega = 1.558039$, while the 5:6 torus whose $W^{s}$ was computed was the one with $\omega = 1.030011$. 

\begin{table} \label{table:gpuBench}
\centering
\begin{tabular}{l | c | c | c}
 & Device 1 &  Device 2 & Device 3  \\
Total program runtime & 25.93 & 10.52 & 12.57  \\
Mean kernel GPU runtime & 0.14278 & 0.02964 & 0.04160  \\
Kernel time \% of total & 61.68\% & 37.07\% &  31.60\%
\end{tabular}
\caption{Benchmarks for GPU-enabled Manifold Intersection Code (all time values in s)}
\end{table}

\begin{table} \label{table:cpuBench}
\centering
\begin{tabular}{l | c | c | c}
 & Device 1 &  Device 2 & Device 3  \\
Total program runtime & 97.12 & 60.96 & 85.90  \\
Mean kernel CPU runtime & 0.79743 & 0.49271 & 0.70800  \\
Kernel time \% of total & 91.96\% & 90.53\% & 92.31 \%
\end{tabular}
\caption{Benchmarks for CPU-only Manifold Intersection Code (all time values in s)}
\end{table}

We timed the Julia program runtime for this application on all three devices. The program carries out checks up to layer 14, with $8$ pairs of half-layers checked per layer, for a total of 112 half-layer pairs checked for intersection during the entire program execution. Each pair of half-layers in turn corresponds to 2,424,307,712 pairs of quads. The resulting program runtimes are given in Table \ref{table:gpuBench}; already, we can see excellent performance on the 2019-era laptop and the older desktop, with the entire manifold meshes and hundreds of billions of quad pairs being checked for intersection in just over 10 seconds. Even device 1, the older laptop, has reasonable performance as well. As described in Section \ref{compImplementSection}, the OpenCL kernel is run once for each pair of half-layers being checked for intersections. Thus, the OpenCL kernel is run 112 times throughout the program execution; we timed all of these kernel runs, and give the mean runtimes in Table \ref{table:gpuBench} as well. 

Although computationally suboptimal, we can force OpenCL.jl to use the CPU cores for multithreaded kernel execution, instead of the GPU. For the sake of comparison, the results of doing so are given in Table \ref{table:cpuBench}. From this, we see that the use of the GPU speeds up kernel execution by a factor of 5.6x for device 1, 16.6x for device 2, and 17x for device 3. For the CPU-only program, the kernel executions make up over 90 percent of the overall program runtime; hence, the kernel speedup achieved through GPU usage results in a large speedup of the program as well. The GPU-enabled program is 3.75x faster on device 1, 5.8x faster on device 2, and 6.8x faster on device 3 than the CPU-only program. 

As a final note, it is instructive to compare the aforementioned results with those of some of our previous work. The first version\cite{kumar2020} of these algorithms for finding manifold intersections was a MATLAB program whose broad phase only consisted of the bounding box test, written as a CPU-only vectorized (and thus parallel) 4D array operation applied to all pairs of quads without any prior grid-based pruning. Running the MATLAB program on device 1, for the same benchmark presented at the beginning of this section, the bounding box test took 8 seconds for each pair of half-layers. The program runtime was thus close to 1000 s. 

The second version\cite{kumar2021feb} was a Julia program using OpenCL.jl and GPUs, which did the bounding box and quick M\"oller tests on the GPU but also did not include the spatial partitioning step in its broad phase. Thus again, the bounding box test was applied to all pairs of quads from each pair of  half-layers. In this second study, we had access to the JPL DGX High Performance Computing platform, from which Julia had could use 16 CPU threads and an Nvidia Tesla V100 GPU with 16GB VRAM; this is far more powerful than any of devices 1, 2, or 3. Nevertheless, the same benchmark from earlier in this section took 16 seconds on DGX, which is worse than both devices 2 and 3! This very clearly illustrates the importance of the spatial partitioning for achieving good algorithm performance. 

\begin{figure}
\includegraphics[width=0.495\columnwidth]{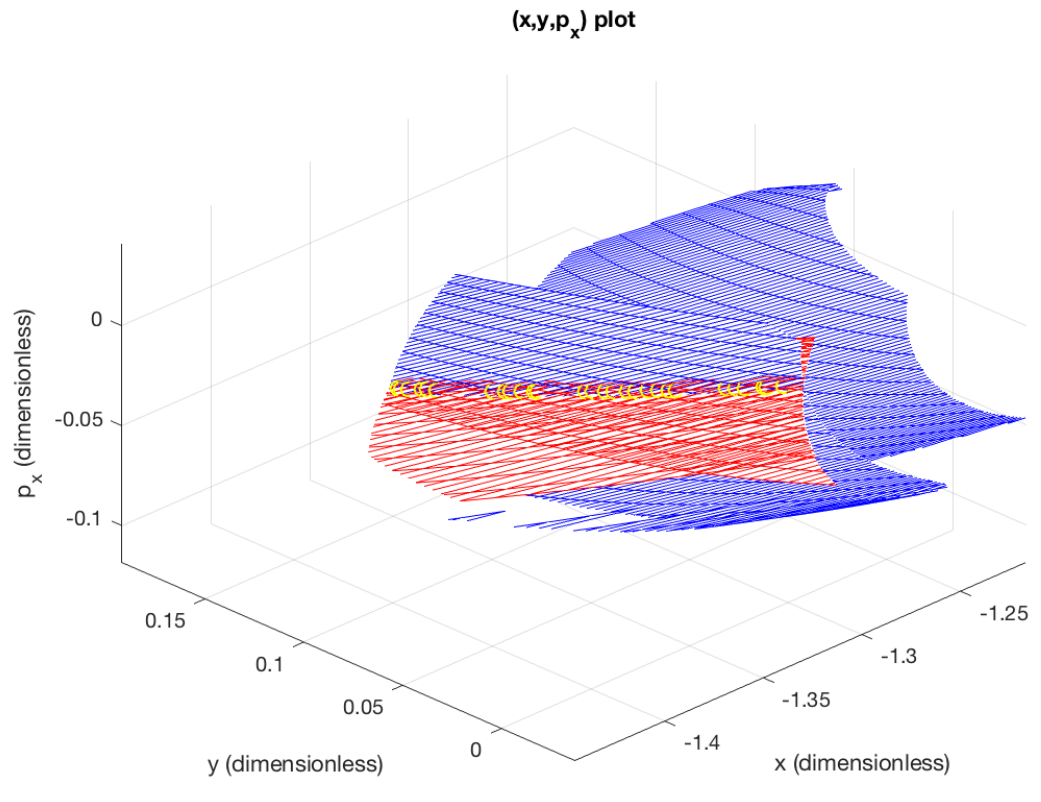}
\includegraphics[width=0.495\columnwidth]{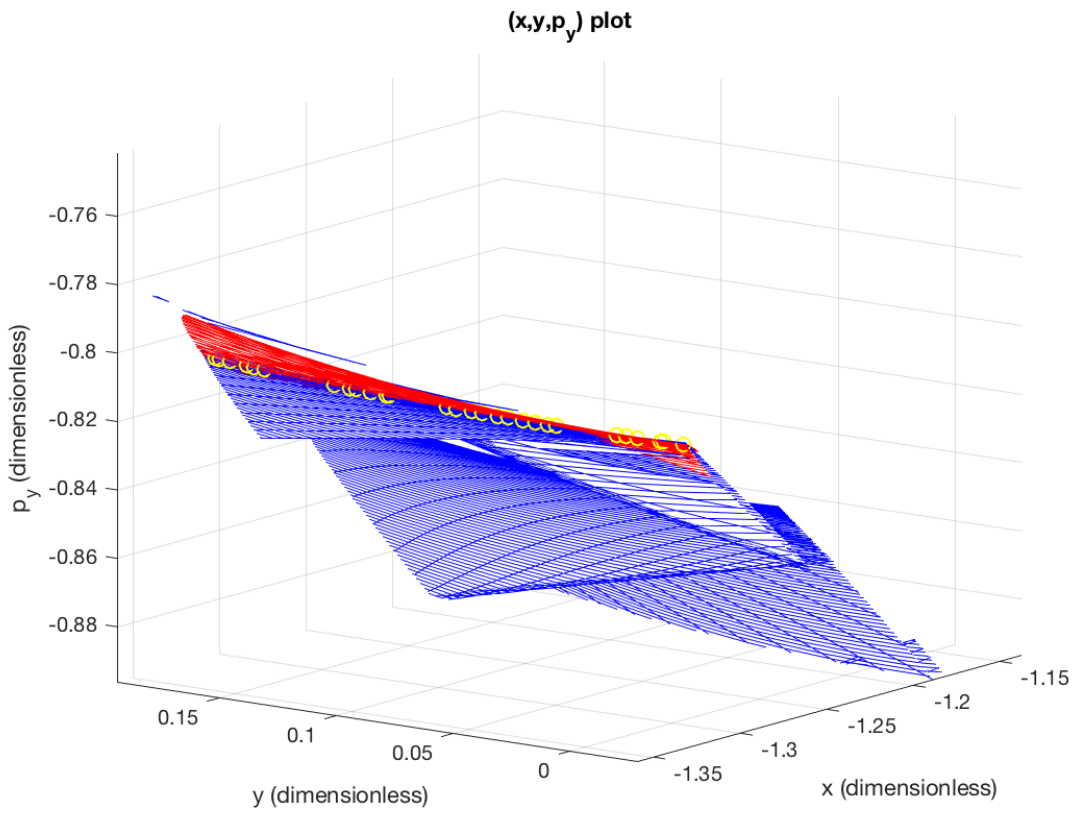}
\caption{ \label{fig:connections} 3:4 $W^{u}$ (red) and 5:6 $W^{s}$ (blue) heteroclinic connections (yellow circles) of resonant tori in Jupiter-Europa PERTBP}
\end{figure}

Finally, we show some manifold mesh intersections output by this test application in Figure \ref{fig:connections}; the plot on the left is the zoomed-in projection onto $(x,y,p_{x})$ space of the intersections, shown as yellow circles. The plot on the right is the projection onto $(x,y,p_{y})$ space of the same intersections. These figures are repeated from our previous study\cite{kumar2020}; since the benchmark (and the numerical results) in this study are the same as the previous ones, newly generated plots of the mesh intersections found using Julia look exactly the same. 

\section{Refinement of Approximate Manifold Intersections} \label{refinementSection}

We have shown how to represent the unstable and stable manifolds $W_{1}^{u}$ and $W_{2}^{s}$ as meshes, and have also given fast methods for finding intersections of these meshes in 4D space. However, these meshes are made of triangles, which linearly interpolate points between their vertices. Of course, this interpolation has error, so an intersection of the meshes is not an exact heteroclinic connection.  We now seek to correct the approximate heteroclinic connections found in the mesh-based search from the previous section. We wish to find solutions $\bold{x} = (\theta_{u}, s_u, \theta_{s}, s_{s})$ of the equation
\begin{equation} \label{tosolve} f(\bold{x}) = f(\theta_{u}, s_u, \theta_{s}, s_{s}) \stackrel{\text{def}}{=}  W^{u}_{1}(\theta_{u}, s_u) -W^{s}_{2}(\theta_{s}, s_{s}) =0\end{equation}
As discussed during the description of the narrow phase in Section \ref{narrowPhase}, we already will have decent initial guesses for $(\theta_{u}, s_u, \theta_{s}, s_{s})$ from the mesh search. Hence, we can use differential correction to solve Equation \eqref{tosolve}, but to do that we must be able to differentiate its LHS. 

Denote $\partial_{\theta} = \frac{\partial }{\partial \theta} $ and $\partial_{s} = \frac{\partial }{\partial s} $. To differentiate the LHS of Equation \eqref{tosolve}, we need the partial derivatives $\partial_{\theta} W^{u}_{1}$, $\partial_{s} W^{u}_{1}$, $\partial_{\theta}  W^{s}_{2}$, and  $\partial_{s} W^{u}_{1}$ evaluated at $(\theta_{u}, s_u, \theta_{s}, s_{s})$. If $s_{u} \notin [-D_{u}, D_{u}]$ or $s_{s} \notin [-D_{s}, D_{s}]$, we cannot just differentiate the Fourier-Taylor parameterizations of the manifolds and evaluate the result at $(\theta_{u}, s_u, \theta_{s}, s_{s})$. However, the parameterizations are still of use. Applying Equation \eqref{globou} to $W^{u}_{1}$ and Equation \eqref{globos} to $W^{s}_{2}$, we have
\begin{equation} \label{wuequation} W^{u}_{1}(\theta, s) = F^{m}(W^{u}_{1}(\theta-m \omega_{1}, \lambda_{u}^{-m}s)) \end{equation}
\begin{equation} \label{wsequation} W^{s}_{2}(\theta, s) = F^{-n}(W^{s}_{2}(\theta+n \omega_{2}, \lambda_{s}^{n}s)) \end{equation}
where $\omega_{1}, \omega_{2}$ are the rotation numbers of $W^{u}_{1}$ and $W^{s}_{2}$. Differentiating equations \eqref{wuequation} and \eqref{wsequation} gives 
\begin{equation} \label{dwudthequation} \partial_{\theta} W_{1}^{u} (\theta, s) = DF^{m}(W_{1}^{u}(\theta-m \omega_{1}, \lambda_{u}^{-m}s)) \, \partial_{\theta} W_{1}^{u} (\theta-m \omega_{1}, \lambda_{u}^{-m}s) \end{equation}
\begin{equation} \label{dwudsequation} \partial_{s} W_{1}^{u} (\theta, s) = \lambda_{u}^{-m} DF^{m}(W_{1}^{u}(\theta-m \omega_{1}, \lambda_{u}^{-m}s))  \,  \partial_{s} W_{1}^{u} (\theta-m \omega_{1}, \lambda_{u}^{-m}s) \end{equation}
\begin{equation} \label{dwsdthequation} \partial_{\theta} W_{2}^{s} (\theta, s) = DF^{-n}(W_{2}^{s}(\theta+n \omega_{2}, \lambda_{s}^{n}s))  \,  \partial_{\theta} W_{2}^{s} (\theta+n \omega_{2}, \lambda_{s}^{n}s) \end{equation}
\begin{equation} \label{dwsdsequation} \partial_{s} W_{2}^{s} (\theta, s) = \lambda_{s}^{n} DF^{-n}(W_{2}^{s}(\theta+n \omega_{2}, \lambda_{s}^{n}s))  \,  \partial_{s} W_{2}^{s} (\theta+n \omega_{2}, \lambda_{s}^{n}s) \end{equation}
Now, if we choose $m$ and $n$ large enough such that $|\lambda_{u}^{-m}s_{u}| < D_{u}$ and $|\lambda_{s}^{n}s_{s}| < D_{s}$, then one can use equations \eqref{dwudthequation}-\eqref{dwsdsequation} to compute the partials at any $(\theta_{u}, s_u, \theta_{s}, s_{s})$. Since $W^{u}_{1}$ has a Fourier-Taylor series parameterization valid for $|s| < D_{u}$, we can directly evaluate $W_{1}^{u}(\theta_{u}-m \omega_{1}, \lambda_{u}^{-m}s_{u}) $. We can also differentiate the parameterization with respect to $\theta$ and $s$ to get Fourier-Taylor series for $\partial_{\theta} W^{u}_{1}$ and $\partial_{s} W^{u}_{1}$, which can then  both be evaluated at $(\theta,s)=(\theta_{u}-m \omega_{1}, \lambda_{u}^{-m}s_{u}) $. Finally, the $DF^{m}$ from equations \eqref{dwudthequation} and \eqref{dwudsequation} is just a state transition matrix,  found by time-$2\pi m/\Omega_{p}$ numerical integration of the variational equations starting at the known point $W^{u}_{1}(\theta_{u}-m \omega_{1}, \lambda_{u}^{-m}s_{u}) $. All this allows us to compute the partials of  $W^{u}_{1}$; the partials of $W^{s}_{2}$ are done very similarly. 

We now describe how to compute the various quantities on the RHS of equations \eqref{dwudthequation}-\eqref{dwsdsequation}. Without loss of generality, we describe the process for $W_{1}^{u}$; $W_{2}^{s}$ is done in the same way. Recall from Equation \eqref{series} that our Fourier-Taylor parameterization $W_{1}^{u}$ is of the form 
\begin{equation} \label{seriesu} W^{u}_{1}(\theta, s) = \sum_{k \geq 0} W^{u}_{1,k}(\theta)s^{k} \end{equation}
 The coefficients $W^{u}_{1,k}(\theta)$ are stored as arrays of their values at $N$ evenly spaced $\theta$ values $\theta_{u,i} = 2\pi i/N$, $i= 0, 1, \dots, N-1$. Hence, given $(\theta_{u}, s_{u})$, we first evaluate $W_{1}^{u}(\theta_{u,i}, \lambda_{u}^{-m}s_{u}) $ at all the $\theta_{u,i}$ using Equation \eqref{seriesu} and the known coefficients $W^{u}_{1,k}(\theta_{u,i})$. Then, trigonometric interpolation\cite{triginterp} allows us to find the value of $W_{1}^{u}(\theta_{u}-m \omega_{1}, \lambda_{u}^{-m}s_{u}) $ needed in Equations \eqref{dwudthequation}-\eqref{dwudsequation}, which is used to start the numerical integration of the state transition matrix $DF^{m}$. 

For the RHS of Equation \eqref{dwudthequation}, $\partial_{\theta}W_{1}^{u}(\theta_{u}-m \omega_{1}, \lambda_{u}^{-m}s_{u}) $ can be found by first using all the $W_{1}^{u}(\theta_{u,i}, \lambda_{u}^{-m}s_{u}) $ values found earlier to compute $\partial_{\theta}W_{1}^{u}(\theta_{u,i}, \lambda_{u}^{-m}s_{u}) $. Trigonometric interpolation then yields the value at $\theta=\theta_{u}-m \omega_{1}$. To get $\partial_{\theta}W_{1}^{u}(\theta_{u,i}, \lambda_{u}^{-m}s_{u}) $ at all the $\theta_{u,i}$ from knowledge of $W_{1}^{u}(\theta_{u,i}, \lambda_{u}^{-m}s_{u}) $ uses a general FFT based technique. In particular, suppose one has a periodic function $a(\theta)$ whose values are stored at $\theta_{i} = 2\pi i/N$, $i= 0, \dots, N-1$. Then, using the usual formula relating $a(\theta_{i})$ to its FFT coefficients $\hat a(i)$, we find that
\begin{equation} \label{translate2} a(\theta_{i})  = \frac{1}{N} \sum_{k =0}^{N-1} \hat a(k) e^{jk\theta_{i}}  \rightarrow \partial_{\theta} a(\theta_{i})  = \frac{1}{N} \sum_{k =0}^{N-1} [jk \hat a(k)] e^{jk\theta_{i}}  \end{equation}
where $j$ denotes $\sqrt{-1}$ in Equation \eqref{translate2}. Hence, given the $a(\theta_{i})$ values, to find the values $\partial_{\theta} a(\theta_{i}) $, one takes the FFT, multiplies the $k$th Fourier coefficient by $jk$, and takes the inverse FFT. This completes the tools required to find $\partial_{\theta} W_{1}^{u} (\theta_{u}-m \omega_{1}, \lambda_{u}^{-m}s_{u})$. 

For the RHS of Equation \eqref{dwudsequation},  one can find $\partial_{s} W_{1}^{u} (\theta_{u}-m \omega_{1}, \lambda_{u}^{-m}s_{u})$ by first differentiating Equation \eqref{seriesu} at each fixed $\theta_{u,i}$ grid value with respect to $s$. This yields $N$ polynomials in $s$
\begin{equation} \label{seriesuds} \partial_{s} W^{u}_{1}(\theta_{u,i}, s) = \sum_{k \geq 0} (k+1)W^{u}_{1,k+1}(\theta_{u,i})s^{k} \end{equation}
with known coefficients. Next, the $\partial_{s} W^{u}_{1}(\theta_{u,i}, s)$ series can be evaluated at $s = \lambda_{u}^{-m}s_{u}$ for all the $\theta_{u,i}$, finally followed by trigonometric interpolation to find $\partial_{s} W_{1}^{u} (\theta_{u}-m \omega_{1}, \lambda_{u}^{-m}s_{u})$. 

With all quantities from the RHS of Equations  \eqref{dwudthequation} and \eqref{dwudsequation} found, we can now compute the desired partials at $(\theta_{u}, s_{u})$. As mentioned earlier, the partials of $W^{s}_{2}$ can be found in the exact same manner, after which we can solve Equation \eqref{tosolve}. With $\bold{x} = (\theta_{u}, s_u, \theta_{s}, s_{s})$ and letting $\bold{x}_{0}$ be the initial guess found earlier for $\bold{x}$ solving Equation \eqref{tosolve}, we use the damped Newton method
\begin{equation} \label{dampedNewtonEquation} \bold{x}_{k+1} = \bold{x}_{k} - \alpha Df^{-1}(\bold{x}_{k}) f(\bold{x}_{k}) \end{equation}
to differentially correct $\bold{x}$ until we have a solution to Equation \eqref{tosolve} within tolerance. Here, $Df = [\partial_{\theta_{u}}f \,\, \partial_{s_{u}}f \,\, \partial_{\theta_{s}}f \,\, \partial_{s_{s}}f ]$, and $0< \alpha < 1$. Trial and error is used to find a value of $\alpha$ such that the iteration converges. We have used $\alpha$ values anywhere from $0.5$ to $0.01$ in our computations. 

\begin{figure}[h]
          \includegraphics[width=0.495\columnwidth]{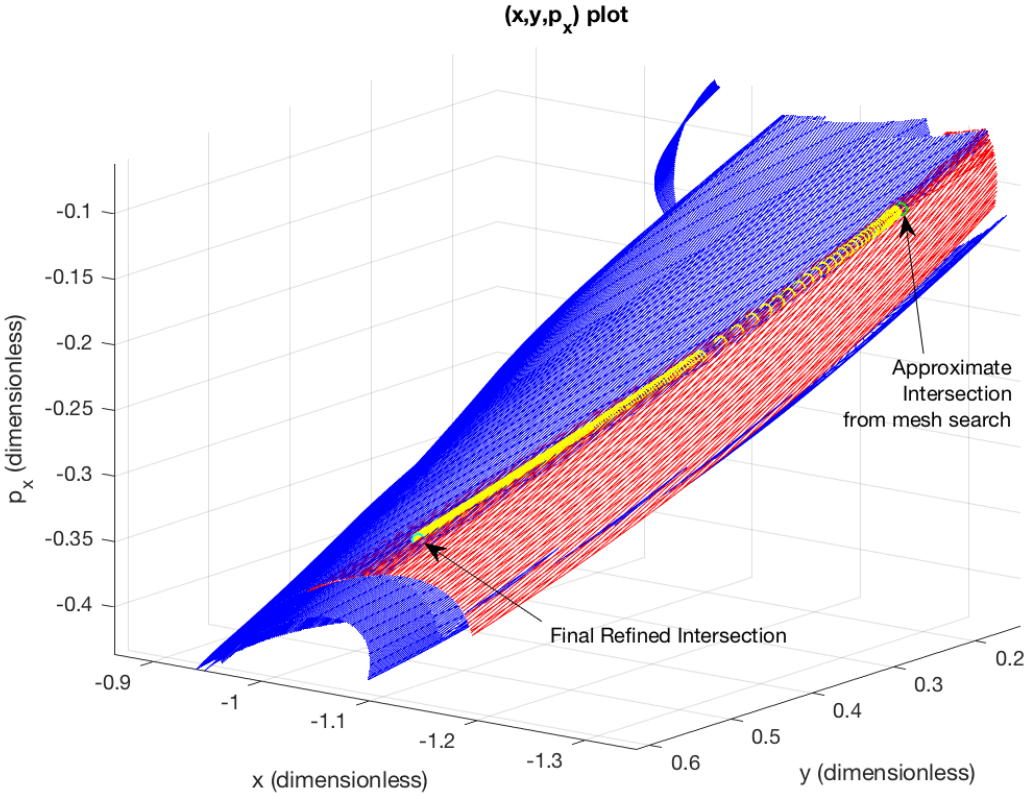}
          \includegraphics[width=0.495\columnwidth]{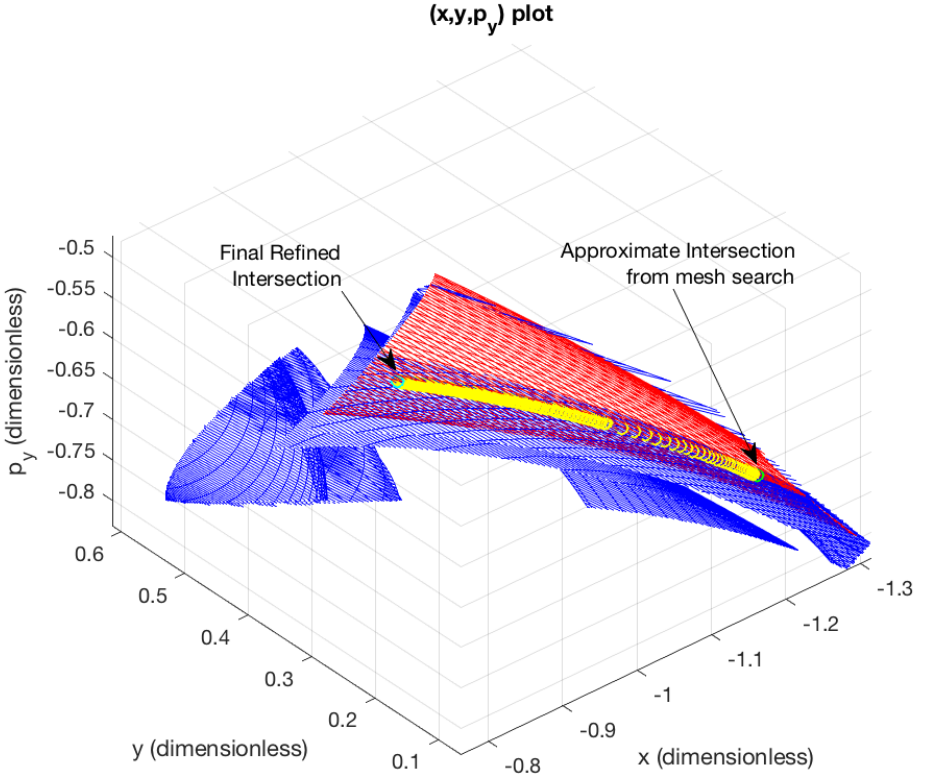}
	          \caption{\label{dampedNewton} 3D projections of iterates of damped Newton method for refinement of an approximate connection found from mesh search (Jupiter-Europa PERTBP 3:4 $W^{u}$ in red, 5:6 $W^{s}$ in blue)}
	 \end{figure}

Using the damped Newton method, we have been able to differentially correct several of the approximate intersections $(\theta_{u}, s_u, \theta_{s}, s_{s})$ found in the mesh search benchmark described in Section \ref{benchmarkSection}, from an error of 0.01 in Equation \eqref{tosolve} to errors of less than $10^{-7}$. An example differential correction is displayed in Figure \ref{dampedNewton}, for one of the connections shown in Figure \ref{fig:connections}. The initial guess in Figure \ref{dampedNewton} is shown in green, the iterates in yellow, and the final converged solution in cyan. We changed the value of $\alpha$ at one point during the iteration, hence the uneven spacing of the iterates. Also, note that the damped Newton iterates move a large distance away from the initial guess found in the mesh search. This is not unexpected; the derivative of the LHS of Equation \eqref{tosolve} is almost singular, since PCRTBP manifold intersections (and hence solutions to Equation \eqref{tosolve} in the PCRTBP case) occur along 1D curves rather than at isolated points. In the PCRTBP case, this implies that the derivative at a solution of the equation would actually be singular. The perturbation in the Jupiter-Europa PERTBP is quite weak, so near-singularity is reasonable to expect. It is because of this that the damped Newton method is necessary.

In addition, we found that all of the manifold intersections shown in Figure \ref{fig:connections}, upon differential correction, actually converged to the same refined solution (the green circle in Figure \ref{dampedNewton})! This is despite all of the solutions from Figure \ref{fig:connections} satisfying Equation \eqref{tosolve} with an error of 0.01 or less. Hence, we see that intersecting the discrete mesh representations of the manifolds may actually find mostly intersections which correspond to near-misses rather than true intersections of the manifolds, especially when the periodic perturbation is weak. This, as well as the distance between the initial guess and the refined solution in Figure \ref{dampedNewton}, demonstrates the importance of carrying out the differential correction in order to find the true intersections. The mesh-based search is necessary in order to quickly narrow down potential points of interest, but the few final exact manifold intersections must be found by solving Equation \eqref{tosolve} directly using high-accuracy methods such as those described here. 

\section{Continuation of Manifold Intersections Through Torus Families}

In the previous section, we showed how to compute a refined heteroclinic connection between two individual  whiskered tori in the phase space of a periodically-forced PCRTBP model's stroboscopic map. However, in such systems, whiskered tori generally appear in 1-parameter families. Indeed, in the unperturbed PCRTBP, unstable periodic orbits occur in 1-parameter families, and most of these periodic orbits persist into the periodically perturbed model as well for sufficiently small perturbation $\varepsilon$. For instance, our previous work \cite{kumar2022} on computing tori in the Jupiter-Europa PERTBP found that the very vast majority of the 3:4 and 5:6 Jupiter-Europa PCRTBP resonant unstable periodic orbits did continue into the PERTBP until the physical value $\varepsilon = 0.0094$ of Europa's eccentricity. Thus, just like the PCRTBP periodic orbits, the corresponding 3:4 and 5:6 PERTBP tori are also present in 1-parameter families; one can take the torus rotation number $\omega$ as the parameter. 

In the PCRTBP, as the autonomous Hamiltonian $H_{0}$ given in Equation \eqref{pcrtbpH} is an integral of motion, any unstable periodic orbit can only have a heteroclinic connection to another unstable periodic orbit which belongs to the same energy level. Thus, as the energy varies along 1-parameter families of periodic orbits, one cannot keep the first periodic orbit fixed and numerically continue the heteroclinic connection along the second orbit family in the PCRTBP, nor vice versa, as the energy would change along the continuation. However, in periodically forced PCRTBP systems, this energy integral no longer exists. Thus, it is possible (and indeed happens, as we show in Section \ref{continuationNumericsSection}) that an unstable torus can have heteroclinic connections to a range of tori in another torus family; hence, between any two 1-parameter families of such tori, there is actually a 2-parameter family of heteroclinic connections. To investigate this family of connections, one could use the tools of Sections \ref{layerSection}-\ref{refinementSection} to compute the heteroclinics from each torus in one family to several tori in the other torus family; however, it is inefficient to repeatedly run the entire program for each pair of tori. Instead, one can use the differential correction methods of Section \ref{refinementSection} to implement a continuation scheme, starting from a connection found between a single pair of tori, one from each family. 

Given a 1-parameter torus family, let $W^{u}_{1f}(\theta_{u}, s_u; \omega_{u})$ be a function which, for each value of $\omega_{u}$, represents the unstable manifold of the torus in that family which has rotation number $\omega_{u}$. In order to signify the dependence of the manifold $W_{1f}^{u}$ on the torus family member from which it emanates, we explicitly include the torus rotation number $\omega_{u}$ as a function parameter. Similarly, let $W^{s}_{2f}(\theta_{s}, s_{s}; \omega_{s})$ represent the stable manifold of the torus of rotation number $\omega_{s}$ belonging to a 1-parameter torus family. The functions $W_{1f}^{u}$ and $W^{s}_{2f}$ can (and must) be defined in a manner ``practically differentiable'' with respect to $\omega_{u}$ and $\omega_{s}$; we will explain what is meant by ``practically differentiable'' in the following Section \ref{WhitneySection}, and describe how to construct such functions $W_{1f}^{u}$ and $W^{s}_{2f}$ in Section \ref{W1fSection}. Now, assume that we already have computed $W_{1f}^{u}$, $W^{s}_{2f}$, and a solution $(\theta_{u}, s_u, \theta_{s}, s_{s})=(\theta_{u,0}, s_{u,0}, \theta_{s,0}, s_{s,0})$ to the equation 
\begin{equation} \label{continuationEquation} W^{u}_{1f}(\theta_{u}, s_u; \omega_{u})=W^{s}_{2f}(\theta_{s}, s_{s}; \omega_{s}) \end{equation}
for parameter values $\omega_{u}=\omega_{u,0}$ and $\omega_{s}=\omega_{s,0}$; this corresponds to a heteroclinic connection between the tori having rotation numbers $\omega_{u,0}$ and $\omega_{s,0}$. We will now seek to find a heteroclinic connection for a slightly different pair of tori belonging to the same two torus families. 

Assume that the new tori have rotation numbers $(\omega_{u,*}, \omega_{s,*}) = (\omega_{u,0}, \omega_{s,0}) + (\Delta \omega_{u}, \Delta \omega_{s})$; we thus want to find a solution to Equation \eqref{continuationEquation} for these new values of the parameters $(\omega_{u}, \omega_{s})$. Given ``practically differentiable'' representations of the functions $W_{1f}^{u}$ and $W^{s}_{2f}$, this is very simple to do. If this is the first continuation step, so that the only solution of Equation \eqref{continuationEquation} already known is the one for $(\omega_{u,0}, \omega_{s,0})$, then one can simply take $(\theta_{u}, s_u, \theta_{s}, s_{s})=(\theta_{u,0}, s_{u,0}, \theta_{s,0}, s_{s,0})$ as an initial guess for the solution. If instead one has already carried out a continuation step, and thus already knows two solutions $(\theta_{u,0}, s_{u,0}, \theta_{s,0}, s_{s,0})$ for $(\omega_{u,0}, \omega_{s,0})$ and $(\theta_{u,1}, s_{u,1}, \theta_{s,1}, s_{s,1})$ for $(\omega_{u,1}, \omega_{s,1})$, then as long as $(\omega_{u,*}, \omega_{s,*})$ lies on the same line as $(\omega_{u,0}, \omega_{s,0})$ and $(\omega_{u,1}, \omega_{s,1})$, a better initial guess can be found by a simple linear predictor
\begin{equation} \label{predictor}
 \begin{bmatrix}
\theta_{u} \\ s_u \\ \theta_{s} \\ s_{s} 
\end{bmatrix} 
= \begin{bmatrix}
\theta_{u,0} \\ s_{u,0} \\ \theta_{s,0} \\ s_{s,0} 
\end{bmatrix}  + \frac{\Delta \omega_{u}}{\omega_{u,0}-\omega_{u,1}} 
\begin{bmatrix}
\theta_{u,0}-\theta_{u,1} \\ s_{u,0}-s_{u,1} \\ \theta_{s,0}-\theta_{s,1} \\ s_{s,0}-s_{s,1}
\end{bmatrix} 
\end{equation} 
A common case where $(\omega_{u,*}, \omega_{s,*})$ will be collinear with $(\omega_{u,0}, \omega_{s,0})$ and $(\omega_{u,1}, \omega_{s,1})$ is when one does continuation by only one of $\omega_{u}$ or $\omega_{s}$, i.e. $\Delta \omega_{s}=0$ or $\Delta \omega_{u}=0$; in the latter case, the quotient $ \frac{\Delta \omega_{u}}{\omega_{u,0}-\omega_{u,1}} $ in Equation \eqref{predictor} should be replaced by $ \frac{\Delta \omega_{s}}{\omega_{s,0}-\omega_{s,1}} $. 

After setting the initial guess, one simply applies the differential correction method of Section \ref{refinementSection} to solve for a more precise solution $(\theta_{u}, s_u, \theta_{s}, s_{s})$ of Equation \eqref{continuationEquation} corresponding to the new parameter values $(\omega_{u,*}, \omega_{s,*})$; the correction should converge as long as $(\Delta \omega_{u}, \Delta \omega_{s})$ were taken sufficiently small. Recall from Section \ref{refinementSection} though that in order to evaluate the manifold function $W_{1f}^{u}$ or $W^{s}_{2f}$ and its derivatives at a given point $(\theta_{u}, s_u)$ or $(\theta_{s}, s_{s})$  for use in differential correction, we need the Fourier-Taylor parameterization of the manifold for $\omega_{u}=\omega_{u,*}$ or $\omega_{s}=\omega_{s,*}$, as well as the value of the multiplier $\lambda_{u}$ or $\lambda_{s}$. Thus, we need a way to find both the stable and unstable manifold Fourier-Taylor series as well as their multipliers for different values of $(\omega_{u}, \omega_{s})$. Furthermore, by definition, continuation requires $W_{1f}^{u}$ and $W^{s}_{2f}$ to be ``practically differentiable'' with respect to $\omega_{u}$ and $\omega_{s}$. In the following Section \ref{WhitneySection}, we discuss the type of ``practical differentiability'' one should expect $W_{1f}^{u}$ and $W^{s}_{2f}$ to have. After that, in Section \ref{W1fSection} we explain how to define the functions $W_{1f}^{u}$ and $W^{s}_{2f}$ in such a manner that the continuity requirement is satisfied, and the Fourier-Taylor series are also easily found. 

\subsection{Whitney Differentiability of Manifold Representations with Respect to $\omega$} \label{WhitneySection}

We now describe in more detail the continuity properties of $W_{1f}^{u}$; those of $W^{s}_{2f}$ will be the same. In the previous discussion, we defined $W^{u}_{1f}(\theta_{u}, s_u; \omega_{u})$ as a function which, for each value of $\omega_{u}$, parameterizes the unstable manifold of the torus having rotation number $\omega_{u}$ in that family. Implicit in this definition, though, is that a torus of rotation number $\omega_{u}$ actually exists in the periodically-perturbed PCRTBP system we are studying. However, from the perturbation theory of Hamiltonian systems \cite{morbyBook}, it is known that tori at rational rotation numbers (resonances) will generically not persist into the perturbed system; in fact, for generic perturbations, there will be an interval of $\omega_{u}$ values around each rational rotation number where tori will not exist. Thus, $W^{u}_{1f}(\theta_{u}, s_u; \omega_{u})$ will be undefined on a dense set of the space of $\omega_{u}$ values. Indeed, given the removal of intervals around each resonance, $W^{u}_{1f}$ will be defined only on a  ``fat Cantor set'' (of positive measure) of $\omega_{u}$ values \cite{poschel1982} where the tori persist. 

Since $W^{u}_{1f}(\theta_{u}, s_u; \omega_{u})$ is not even defined on any open interval of $\omega_{u}$ values, its regularity with respect to $\omega_{u}$ is subtle to define. Since the set of allowed $\omega_{u}$ values is a metric space, it is standard to define continuous or even Lipschitz functions, but higher regularity notions are not straightforward. These higher differentiability properties are useful for numerics, though, since they allow to estimate the error in extrapolations and interpolations. The appropriate concept of smoothness to use in sets that do not include balls is that of Whitney differentiability\cite{whitneyExtension,Stein70}. The Whitney extension theorem shows that under certain conditions, a function with prescribed values on a closed set of Euclidean space can be extended to the full Euclidean space in a smooth manner; a function satisfying these conditions is called Whitney differentiable (or smooth). From the practical point of view, this means that for Whitney smooth functions, one can use the usual methods of interpolation and extrapolation. 

As described earlier, we are able to prescribe values of $W^{u}_{1f}(\theta_{u}, s_u; \omega_{u})$ on a (closed) Cantor set of $\omega_{u}$ values; it has been shown that this can be done in a Whitney differentiable manner\cite{poschel1982, chierchiaGallavotti, broerHuitemaSevryuk}, so the resulting $W^{u}_{1f}$ is Whitney smooth in $\omega_u$ for each fixed $\theta_u, s_u$. It is also smooth in the usual sense in $\theta_u, s_{u}$ for each fixed $\omega_u$; it thus follows\cite{Stein70, Llave92} that $W^{u}_{1f}$ is also jointly Whitney smooth in the three variables $(\theta_u, s_u, \omega_u)$ and, therefore, extended smoothly in all the variables. The operator from the original function to its Whitney extension can also be made linear\cite{Stein70}. From the practical point of view, we can consider the function $W^u_{1f}$ as smooth and thus use interpolation to find $W^{u}_{1f}$ at values for which the torus can be computed. 

This Whitney differentiability not only extends the values of our function $W^{u}_{1f}(\theta_{u}, s_u; \omega_{u})$ in a smooth manner, but also justifies the heteroclinic connection continuation procedure by enabling the use of the implicit function theorem to define $(\theta_{u}, s_u, \theta_{s}, s_{s})$ as functions of $(\omega_{u}, \omega_{s})$. In addition, thanks to the Whitney differentiability of the function $W^{u}_{1f}(\theta_{u}, s_u; \omega_{u})$ with respect to $\omega_{u}$, we can expect $W^{u}_{1f}$ to be locally well-approximated by polynomials in $\omega_{u}$. Thus, if we compute $W^{u}_{1f}(\theta_{u}, s_u; \omega_{u})$ on a grid of $\omega_{u}$ values $\{\omega_{u,r}\}_{r=1}^{N_{\omega}}$ it should be possible to interpolate values of $W^{u}_{1f}$ quite accurately between the grid values. We will use this approach in the following section when showing how to practically construct and compute $W^{u}_{1f}$. 

\subsection{Constructing Parameter-Dependent Manifold Representations} \label{W1fSection}

In this section, we will describe how to construct the function $W^{u}_{1f}(\theta_{u}, s_u; \omega_{u})$; the construction of $W^{s}_{2f}(\theta_{s}, s_s; \omega_{s})$ is done in essentially the same manner. Given a family of tori having rotation numbers $\omega_{u} \in [\omega_{u,min}, \omega_{u,max}]$, for each fixed value of $\omega_{u}$, we want $W^{u}_{1f}(\theta_{u}, s_u; \omega_{u})$ to be a solution of Equation \eqref{invariancequationfinal} mapping $\mathbb{T} \times \mathbb{R}$ onto the unstable manifold of the torus of rotation number $\omega_{u}$ belonging to that family. Each torus' unstable multiplier $\lambda_{u}$ will also be dependent on $\omega_{u}$, so we write $\lambda_{u} = \lambda_{u}(\omega_{u})$. In practice, we are able to compute $\lambda_{u}$ values and Fourier-Taylor parameterizations of the torus unstable manifolds on a discrete set of $N_{\omega}$ rotation numbers $\{\omega_{u,r}\}_{r=1}^{N_{\omega}} \subset [\omega_{u,min}, \omega_{u,max}]$. We would now like to use these parameterizations and $\lambda_{u}$ values alongside Equation \eqref{wuequation} to define the values of $W^{u}_{1f}(\theta_{u}, s_u; \omega_{u})$ for other $\omega_{u} \in  [\omega_{u,min}, \omega_{u,max}]$.

To facilitate the computation of the required partial derivatives for differential correction, we will express $W^{u}_{1f}(\theta_{u}, s_u; \omega_{u})$ near $s_{u}=0$ as a parameter-dependent Fourier-Taylor series
\begin{equation} \label{parameterSeries} W^{u}_{1f}(\theta_{u}, s_{u}; \omega_{u}) =\sum_{k \geq 0} W^{u}_{1f,k}(\theta_{u}; \omega_{u}) s_{u}^{k} \end{equation}
For larger $s_{u}$, the values of $W^{u}_{1f}(\theta_{u}, s_u; \omega_{u})$ can then be found by substituting this Fourier-Taylor series and $\lambda_{u}(\omega_{u})$ into Equation \eqref{wuequation}.
Thus, we need to compute $\lambda_{u}$ and the coefficients $W^{u}_{1f,k}(\theta_{u}; \omega_{u})$. A natural way of practically estimating their values for $\omega_{u} \in [\omega_{u,min}, \omega_{u,max}]$ would be to set the values of $W^{u}_{1f,k}(\theta_{u}; \omega_{u,r})$ and $\lambda_{u}(\omega_{u,r})$ using the Fourier-Taylor coefficients and $\lambda_{u}$ already known at $\omega_{u} \in \{\omega_{u,r}\}_{r=1}^{N_{\omega}}$, followed by interpolation for intermediate $\omega_{u}$ values. This works well for $\lambda_{u}(\omega_{u})$. However, there are two underdeterminacies in the computation of torus unstable manifold Fourier-Taylor parameterizations which can prevent $W^{u}_{1f,k}(\theta_{u}; \omega_{u})$ from being Whitney with respect to $\omega_{u}$, and thus cause problems in both interpolation and continuation, unless they are properly handled when setting the values of $W^{u}_{1f,k}(\theta_{u}; \omega_{u,r})$. 

The first underdeterminacy is that given a parameterization $W(\theta,s)$ satisfying Equation \eqref{invariancequationfinal}, $W_{\rho}(\theta,s)=W(\theta+\rho,s)$ will also be a solution for any $\rho \in \mathbb{T}$. Both solutions will map $\mathbb{T} \times \mathbb{R}$ onto the same manifold in $\mathbb{R}^{4}$, but for a fixed $\theta_{0} \in \mathbb{T}$, the points $W(\theta_{0},s)$ and $W_{\rho}(\theta_{0},s)$ may be far apart. The implication of this is that if we compute functions $W^{u}_{r}(\theta,s)$ and $W^{u}_{r+1}(\theta,s)$ representing the unstable manifolds of tori at two rotation numbers $\omega_{u,r}$ and $\omega_{u,r+1}$ respectively, then even for $|\omega_{u,r}-\omega_{u,r+1}|$ arbitrarily small there is no guarantee that $\| W^{u}_{r}(\theta, s)-W^{u}_{r+1}(\theta, s) \|$ will be small. Thus, one cannot simply set $W^{u}_{1f}(\theta_{u}, s_u; \omega_{u,r})=W^{u}_{r}(\theta_{u},s_{u})$ and $W^{u}_{1f}(\theta_{u}, s_u; \omega_{u,r+1})=W^{u}_{r+1}(\theta_{u},s_{u})$ and expect the resulting $W^{u}_{1f}$ to be Whitney in $\omega_{u}$. Instead, given the sequence $\{\omega_{u,r}\}_{r=1}^{N_{\omega}}$ of rotation numbers for which we have computed corresponding manifold parameterizations $W^{u}_{r}(\theta,s) = \sum_{k \geq 0} W^{u}_{r,k}(\theta)s^{k}$, one should:

\begin{enumerate}
\item Set $W^{u}_{1f}(\theta_{u}, s_u; \omega_{u,1})=W^{u}_{1}(\theta_{u},s_{u})$ by setting $W^{u}_{1f,k}(\theta_{u}; \omega_{u,1})=W^{u}_{1,k}(\theta_{u})$ for all $k \geq 0$
\item For $r = 2, \dots, N_{\omega}$, recursively set $W^{u}_{1f}(\theta_{u}, s_u; \omega_{u,r})=W^{u}_{r}(\theta_{u} + \rho_{r},s_{u})$ by setting $W^{u}_{1f,k}(\theta_{u}; \omega_{u,r})=W^{u}_{r,k}(\theta_{u}+\rho_{r})$ for all $k \geq 0$, where the phase shift $\rho_{r}$ is given by
\begin{equation} \label{rhoEquation} \rho_{r} = \argmin_{\alpha \in \mathbb{T}} \int_{\mathbb{T}}\| W^{u}_{r,0}(\theta_{u} + \alpha) - W^{u}_{1f,0}(\theta_{u}; \omega_{u,r-1}) \| \, d\theta_{u} \end{equation}
\end{enumerate}

As mentioned in Section \ref{quasiNewton}, $W^{u}_{r,0}(\theta_{u})$ simply corresponds to the base invariant torus from which the unstable manifold parameterized by $W^{u}_{r}(\theta_{u},s_{u})$ emanates;  thus, what the previous procedure does is to shift the angular phasing of the functions $W^{u}_{r}$ parameterizing the manifolds of the $\omega_{u,r}$ tori so that they all line up with the $\theta_{u}$ phasing of the $\omega_{u,1}$ torus' parameterization. Also, recall that on the computer, functions of $\theta_{u}$ are discretized on a grid of $N$ values $\theta_{u,i} = 2\pi i/N$, $i= 0, 1, \dots, N-1$. Hence, we can approximate Equation \eqref{rhoEquation} as
\begin{equation} \label{rhoEquationDiscrete} \rho_{r} = \argmin_{\alpha \in \left\{ \frac{2\pi i}{N} \right\}_{i=0}^{N-1} } \sum_{i=0}^{N-1} \| W^{u}_{r,0}(\theta_{u,i} + \alpha) - W^{u}_{1f,0}(\theta_{u,i}; \omega_{u,r-1}) \| \end{equation}
since this sum is much easier to compute. Restricting the argument $\alpha$ in the $\argmin$ to the set $\left\{ \frac{2\pi i}{N} \right\}_{i=0}^{N-1}$, which is the same grid as used for $\theta_{i}$, allows us to use the previously computed points for $W^{u}_{r,0}(\theta_{u,i})$ when finding the values of $W^{u}_{r,0}(\theta_{u,i} + \alpha)$ without trigonometric interpolation. We thus used Equation \eqref{rhoEquationDiscrete} in our implementation of this phase shift. 

Once the underdeterminacy in the phase of $\theta_{u}$ has been addressed using the aforementioned method, the $s_{u}^{0}$ coefficient $W^{u}_{1f,0}(\theta_{u}; \omega_{u})$ of the Fourier-Taylor series for $W^{u}_{1f}$ will be Whitney in $\omega_{u}$. However, there is also an underdeterminacy in the scaling of $s_{u}$ which, if not handled, can still prevent the $s^{1}$ (and higher order) Fourier-Taylor coefficients from being Whitney with respect to $\omega_{u}$. In particular, if $W(\theta,s)$ is a parameterization satisfying Equation \eqref{invariancequationfinal}, $W_{L}(\theta,s)=W(\theta,Ls)$ will also be a solution for any $L \in \mathbb{R}$. Again, both solutions will have the same image in $\mathbb{R}^{4}$, but for a fixed $s_{0} \in \mathbb{T}$, the points $W(\theta,s_{0})$ and $W_{L}(\theta,s_{0})$ may be far apart; hence, using the same notation as earlier, even for $|\omega_{u,r}-\omega_{u,r+1}|$ arbitrarily small there is no guarantee that $\| W^{u}_{1f}(\theta, s; \omega_{u,r})-W^{u}_{1f}(\theta, s; \omega_{r+1}) \|$ will be small for $s \neq 0$ (the $s=0$ discontinuity was eliminated by the phase adjustment procedure of the last paragraph).  

Recall the Fourier-Taylor series for $W^{u}_{1f}(\theta_{u}, s_{u}; \omega_{u})$ given in Equation \eqref{parameterSeries}, whose coefficients' values were set earlier on the $\omega_{u}$ grid $\{\omega_{u,r}\}_{r=1}^{N_{\omega}}$. To resolve the discontinuity caused by underdeterminacy in the scaling of $s_{u}$, one should simply replace each of the $N_{\omega}$ preset Fourier-Taylor parameterizations $W^{u}_{1f}(\theta_{u}, s_{u}; \omega_{u,r})$ with a rescaled version $W^{u}_{1f}(\theta_{u}, L_{r}s_{u}; \omega_{u,r})$, where $L_{r}=\| W^{u}_{1f,1}(0; \omega_{u,r}) \|^{-1} $. This is equivalent to replacing each coefficient $W^{u}_{1f,k}(\theta_{u}; \omega_{u,r})$ with $L_{r}^{k}W^{u}_{1f,k}(\theta_{u}; \omega_{u,r})$ for all $k \geq 0$; the new $s_{u}^{1}$ coefficient $L_{r}^{1}W^{u}_{1f,1}(\theta_{u}; \omega_{u,r})= \frac{W^{u}_{1f,1}(\theta_{u}; \omega_{u,r})}{ \| W^{u}_{1f,1}(0; \omega_{u,r}) \| }$ will have unit norm at $\theta_{u}=0$ for all $r=1, \dots, N_{\omega}$. Making an abuse of notation, we henceforth refer to the rescaled Fourier-Taylor series as $W^{u}_{1f}(\theta_{u}, s_{u}; \omega_{u}) = \sum_{k \geq 0} W^{u}_{1f,k}(\theta_{u}; \omega_{u}) s_{u}^{k}$. What the rescaling ensures is that the norm of $W^{u}_{1f,1}(\theta_{u}; \omega_{u})$ will vary continuously with $\omega_{u}$ at $\theta_{u}=0$. Since the norm of $W^{u}_{1f,1}(\theta_{u}; \omega_{u})$ at any one $\theta_{u}$ value determines\cite{kumar2022} the scaling of all the $W^{u}_{1f,k}(\theta_{u}; \omega_{u})$ for all $\theta_{u} \in \mathbb{T}$, the new $W^{u}_{1f,k}(\theta_{u}; \omega_{u})$ will also all be continuously scaled in $\omega_{u}$. In our experience, the resulting $W^{u}_{1f}(\theta_{u}, s_{u}; \omega_{u})$ will be Whitney in $\omega_{u}$ as desired.

\subsection{Numerical Demonstration and Results} \label{continuationNumericsSection}

We were able to apply the continuation methods just described on a heteroclinic connection between the same pair of Jupiter-Europa PERTBP 3:4 and 5:6 resonant tori used for numerical demonstrations in the previous Sections \ref{benchmarkSection} and \ref{refinementSection}. We started from an intersection between their manifolds at the phase space point $(x,y,p_{x},p_{y}) = (-0.96064, 0.88783, -0.51377, -0.64714)$, which corresponds to the solution $(\theta_{u}, s_u, \theta_{s}, s_{s})= (2.39703, -77.73428, 1.83093,	202.62277)$ of Equation \eqref{continuationEquation} with parameter values $\omega_{u} = 1.558039$ and $\omega_{s} = 1.030011$. The parameterizations were scaled in $s_{u}$ and $s_{s}$ as described at the end of the previous section, and their $\theta_{u}=0$ and $\theta_{s}=0$ points were taken to lie on the negative x-axis (which fixes the rotational phasing of the parameterizations). 

\begin{figure}[h]
	\begin{centering}
          \includegraphics[width=0.495\columnwidth]{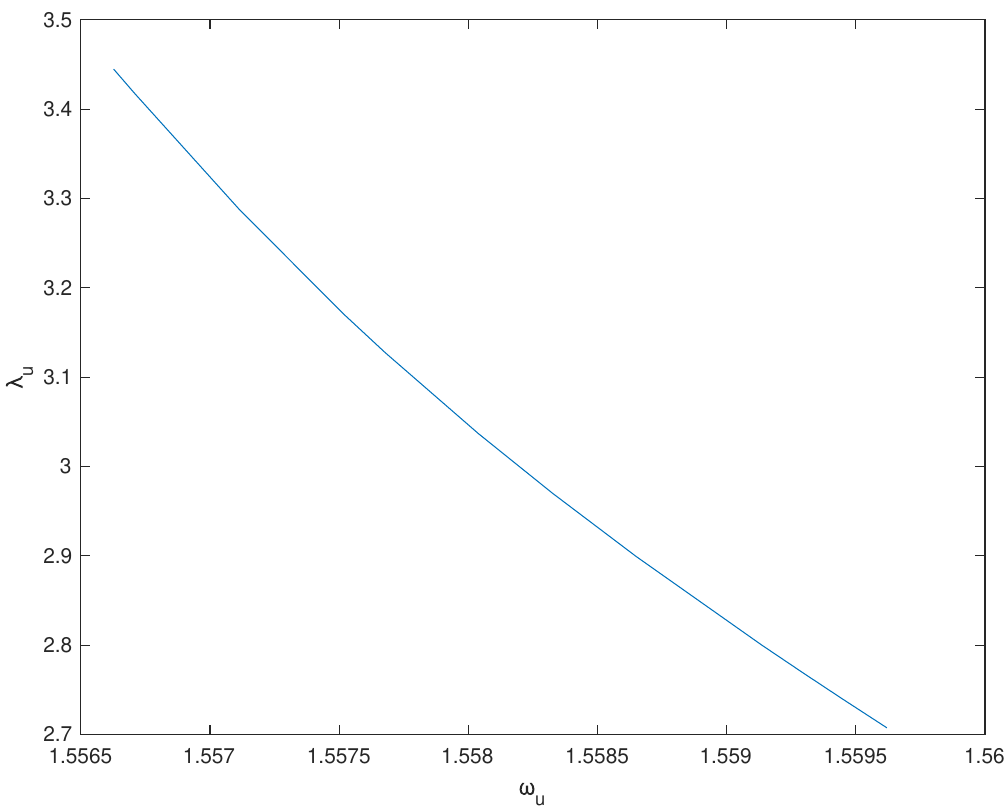}
	          \caption{\label{lambdaVsOmega} $\lambda_{u}$ vs $\omega_{u}$ for the 3:4 Jupiter-Europa PERTBP torus family}
	  \end{centering}
	 \end{figure}

From this starting point, we did a continuation by changing $\omega_{u}$ while keeping $\omega_{s}$ fixed. Thus, we first had to construct $W^{u}_{1f}(\theta_{u}, s_{u}; \omega_{u})$ and $\lambda_{u}(\omega_{u})$ as described in Section \ref{W1fSection}. We computed the Fourier-Taylor parameterizations and $\lambda_{u}$ values for several tori with rotation numbers in the range $\omega_{u} \in [1.55663, 1.55962]$. Plotting the curve of the computed values of $\lambda_{u}$ versus $\omega_{u}$ in Figure \ref{lambdaVsOmega}, we see that the resulting function should indeed be well approximated by a polynomial interpolation; we used quadratic interpolation in this study. 

\begin{figure}[h]
          \includegraphics[width=0.495\columnwidth]{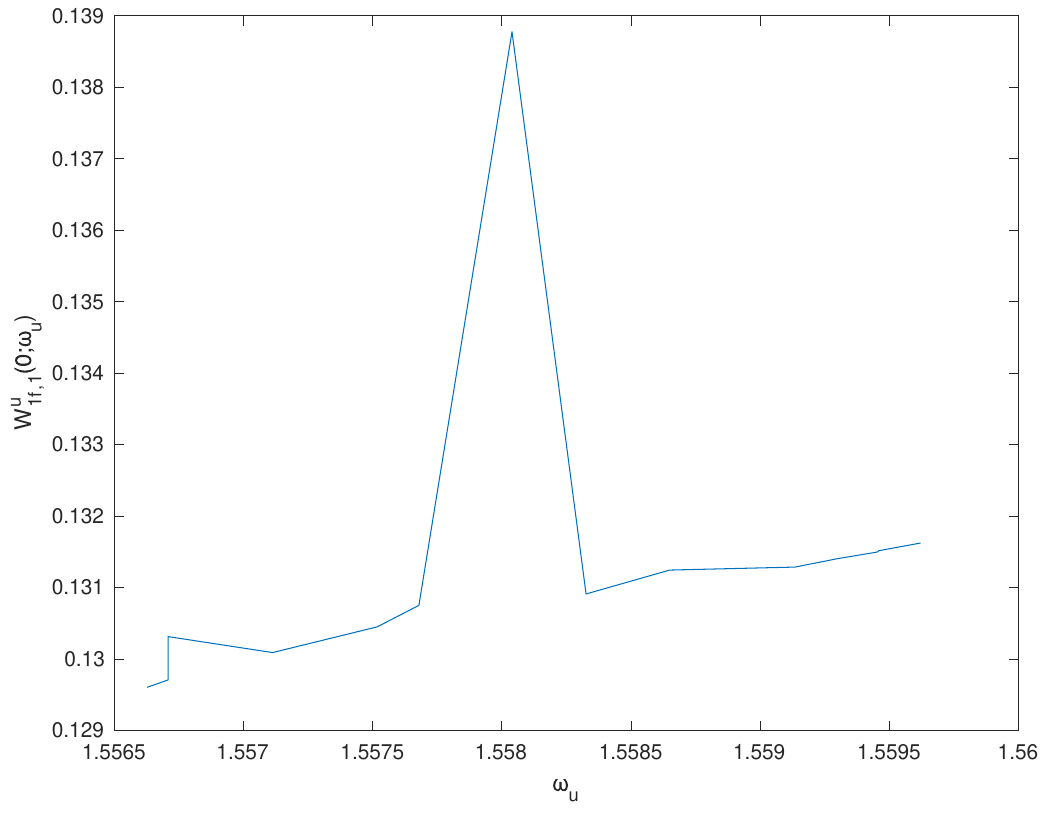}
                    \includegraphics[width=0.495\columnwidth]{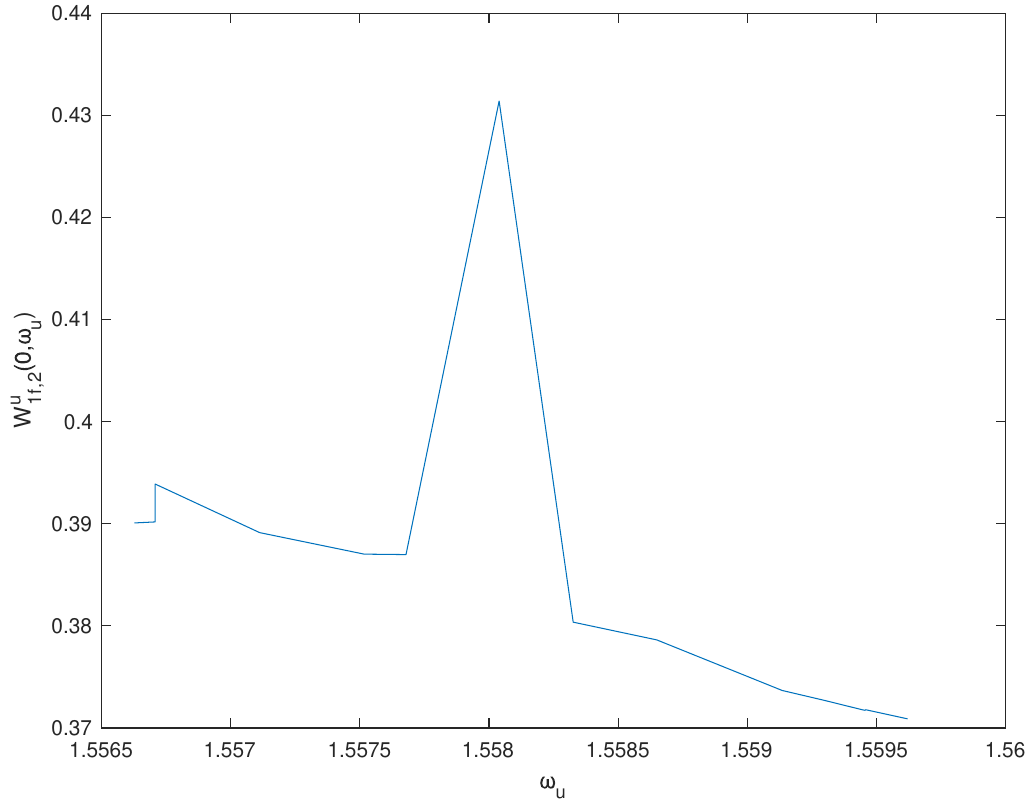}
	          \caption{\label{coeffsVsOmegaUnscaled} $x$-component of unscaled $W^{u}_{1f,1}(0; \omega_{u})$ (left) and $W^{u}_{1f,2}(0; \omega_{u})$ (right) vs $\omega_{u}$}
	 \end{figure}

\begin{figure}[h]
          \includegraphics[width=0.495\columnwidth]{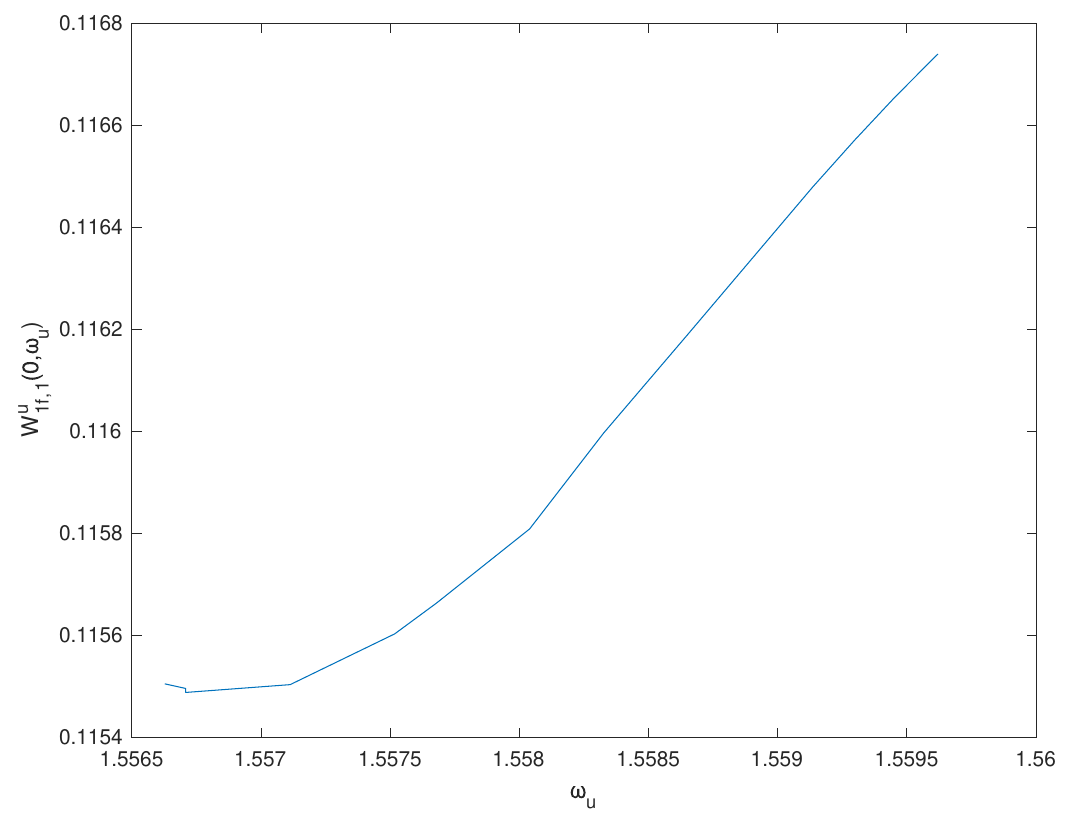}
                    \includegraphics[width=0.495\columnwidth]{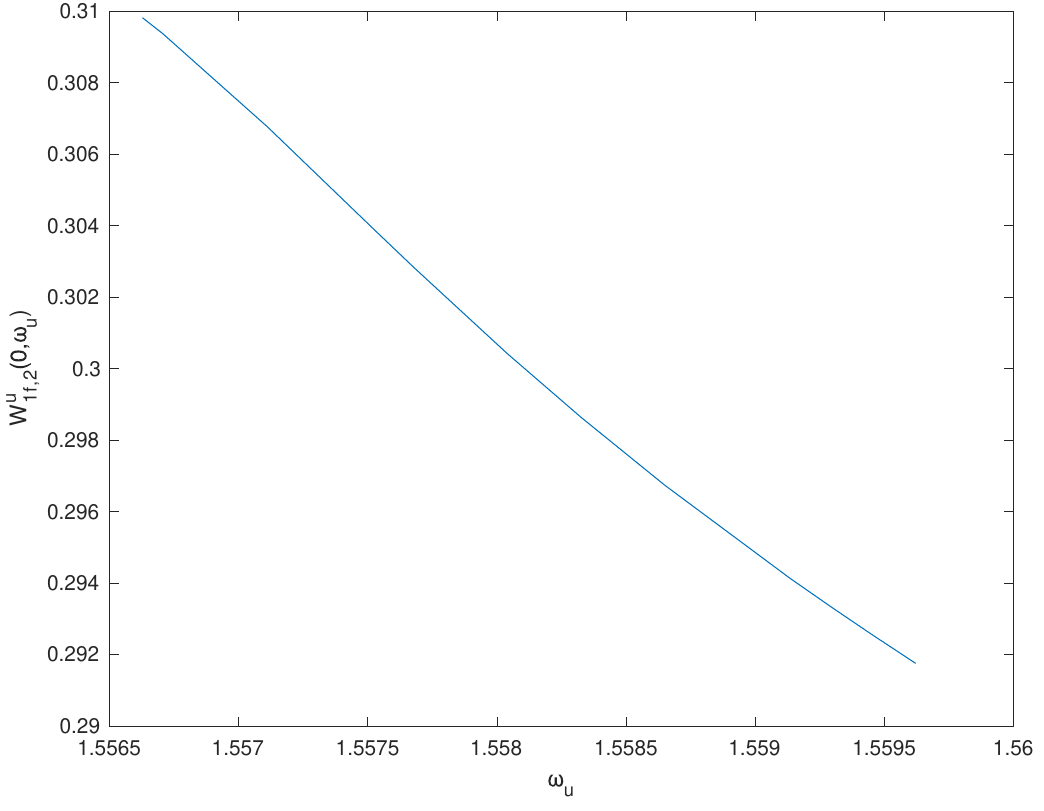}
	          \caption{\label{coeffsVsOmegaScaled} $x$-component of rescaled $W^{u}_{1f,1}(0; \omega_{u})$ vs $\omega_{u}$ (left) and $W^{u}_{1f,2}(0; \omega_{u})$ (right)}
	 \end{figure}

Trying to plot the Fourier-Taylor coefficients versus $\omega_{u}$ after the phase shift, but without the $s_{u}$ rescaling, however, yields a seemingly discontinuous curve, which we show in Figure \ref{coeffsVsOmegaUnscaled}; there, we plot the $x$-component at $\theta_{u}=0$ of the $s_{u}^{1}$ coefficient $W^{u}_{1f,1}(\theta_{u}; \omega_{u})$ and the $s_{u}^{2}$ coefficient $W^{u}_{1f,1}(\theta_{u}; \omega_{u})$ versus $\omega_{u}$. After the rescaling, however, the plots of the same coefficients versus $\omega_{u}$, shown in Figure \ref{coeffsVsOmegaScaled}, are far better behaved, and are clearly amenable to polynomial interpolation just like $\lambda_{u}$. 
	 As $\omega_{s}$ is not being changed, one does not need to construct $W^{s}_{2f}(\theta_{s}, s_{s}; \omega_{s})$ or $\lambda_{s}(\omega_{s})$ apart from the single manifold computation at $\omega_{s} = 1.030011$, which was already done in preparation for the benchmark in Section \ref{benchmarkSection}. 
	 
\begin{figure}[h]
          \includegraphics[width=0.495\columnwidth]{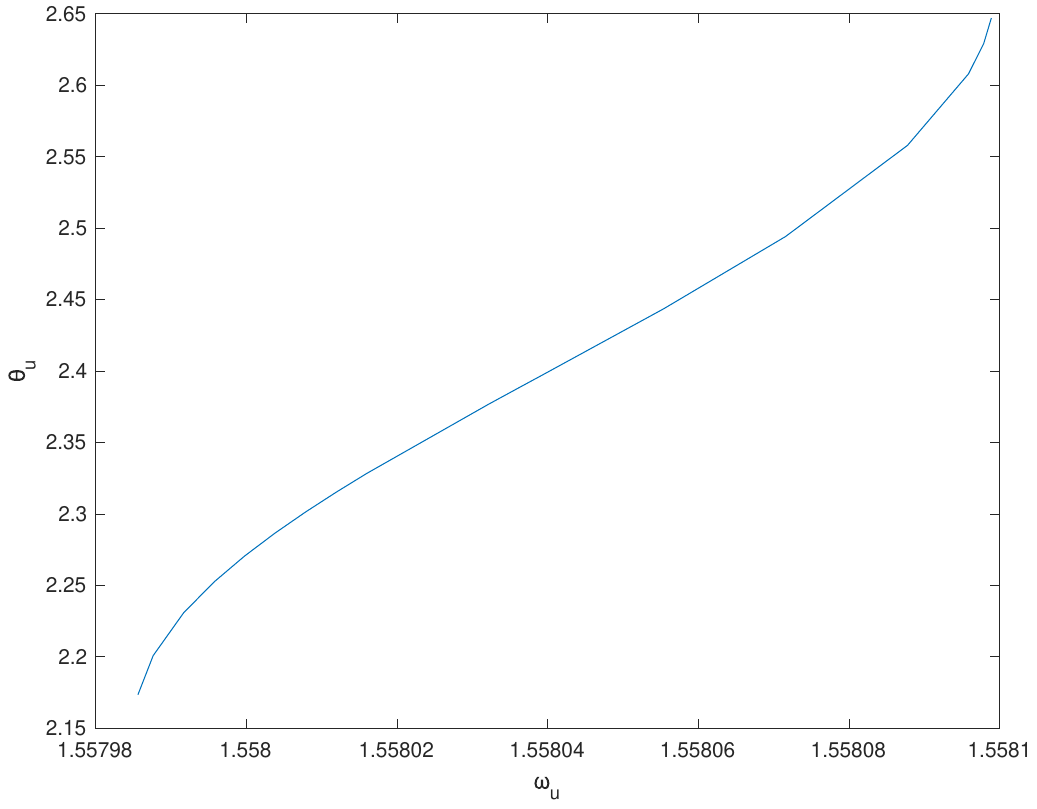}
                    \includegraphics[width=0.495\columnwidth]{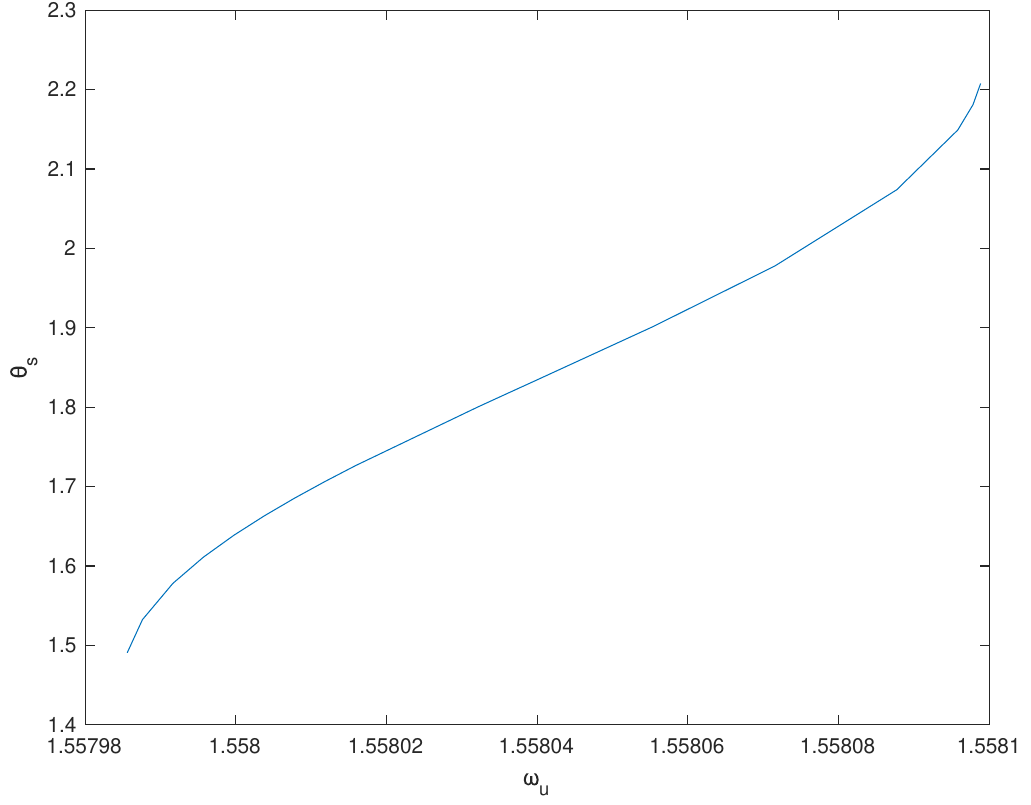}
	          \caption{\label{solnsThetaVsOmega} $\theta_{u}$ (left) and $\theta_{s}$ (right) of solutions of Equation \eqref{continuationEquation} vs $\omega_{u}$}
	 \end{figure}

With the functions $W^{u}_{1f}(\theta_{u}, s_{u}; \omega_{u})$ and $\lambda_{u}(\omega_{u})$ constructed, the continuation by $\omega_{u}$ was then carried out. We successfully computed a solution curve in $(\theta_{u}, s_u, \theta_{s}, s_{s})$ space for $\omega_{u} \in [1.55799, 1.55810]$; the $\theta_{u}$ and $\theta_{s}$ components of this curve are plotted in Figure \ref{solnsThetaVsOmega}, while its $s_{u}$ and $s_{s}$ components are plotted in Figure \ref{solnsSvsOmega}. 
	\begin{figure}[h]
          \includegraphics[width=0.495\columnwidth]{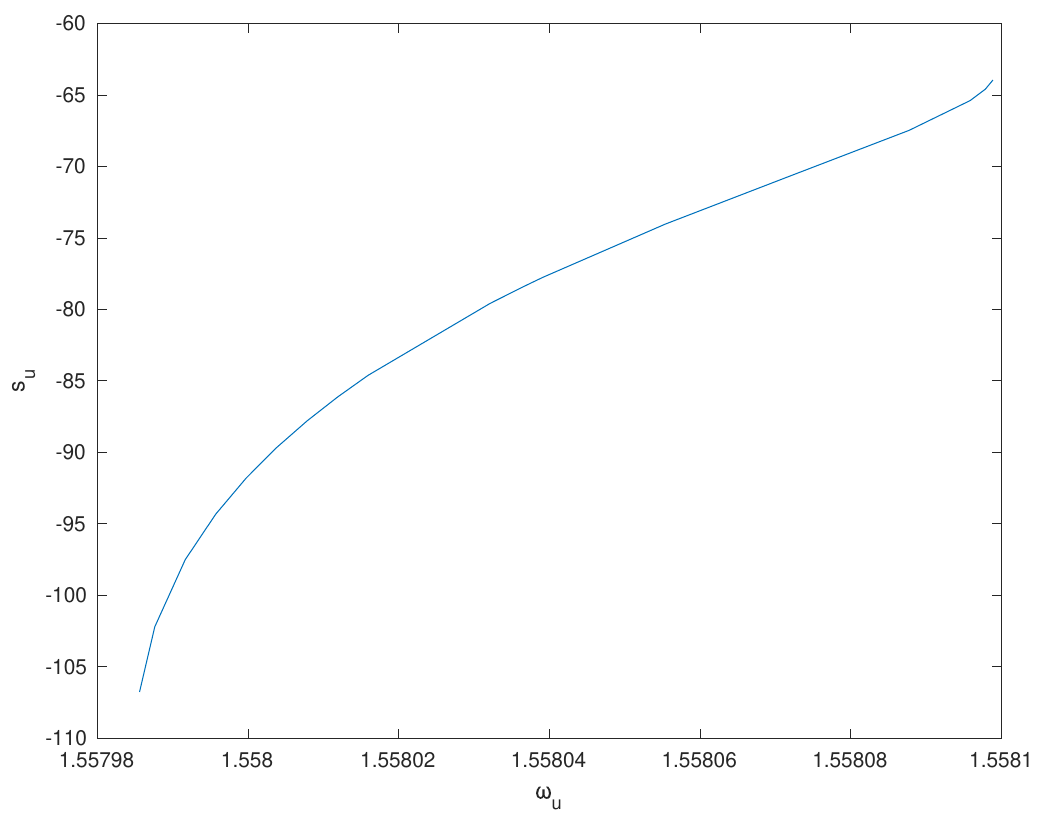}
                    \includegraphics[width=0.495\columnwidth]{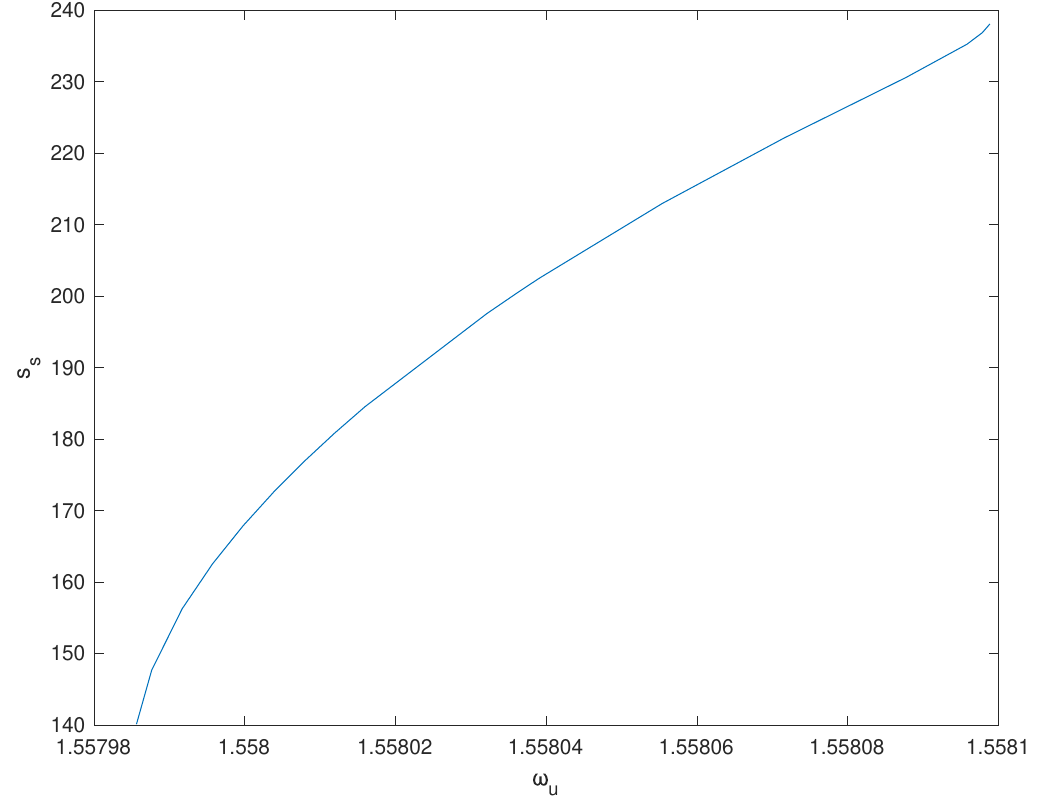}
	          \caption{\label{solnsSvsOmega} $s_{u}$ (left) and $s_{s}$ (right) of solutions of Equation \eqref{continuationEquation} vs $\omega_{u}$}
	 \end{figure}
We can see that these curves are also well behaved and amenable to interpolation, which helps justify the linear predictor of Equation \eqref{predictor} for setting initial guesses to help solve Equation \eqref{continuationEquation}. 3D projections of the heteroclinic intersection points computed using the continuation method are plotted in $(x,y,p_{x})$ and $(x,y,p_{y})$ Cartesian space in Figure \ref{continuationCartesian}; in that figure we also plot the destination 5:6 stable manifold which contains all the computed points (as the continuation involved changing $\omega_{u}$, no single 3:4 unstable manifold contains all of these points, so we do not display any 3:4 manifold). Note that even for this small range of rotation numbers $\omega_{u} \in [1.55799, 1.55810]$ of the origin 3:4 tori, the intersections vary greatly in their Cartesian positions.
	\begin{figure}[h]
          \includegraphics[width=0.495\columnwidth]{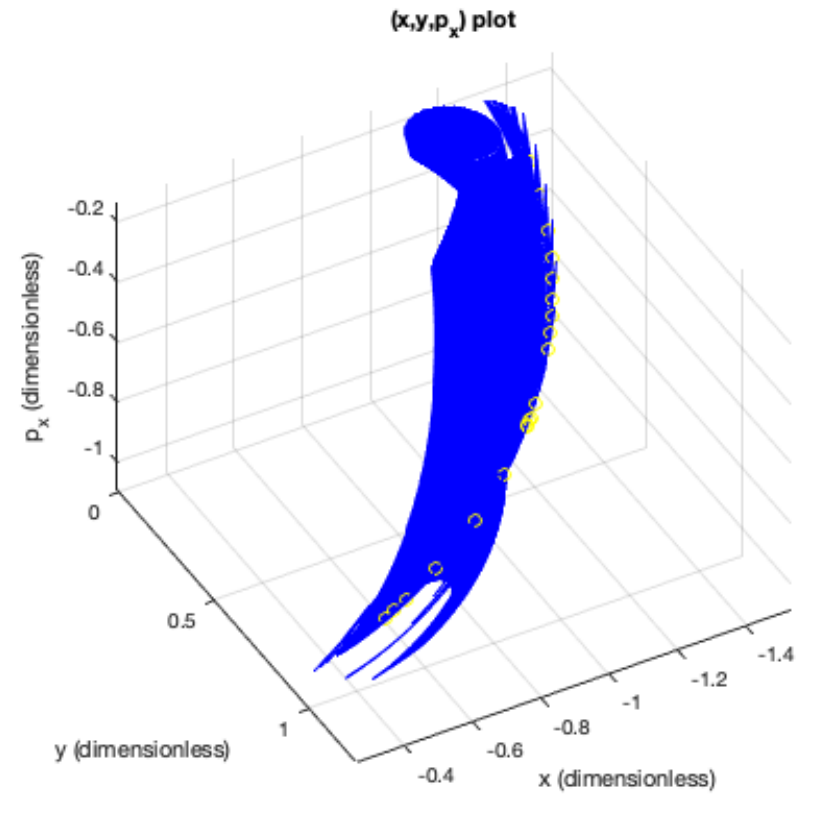}
          \includegraphics[width=0.495\columnwidth]{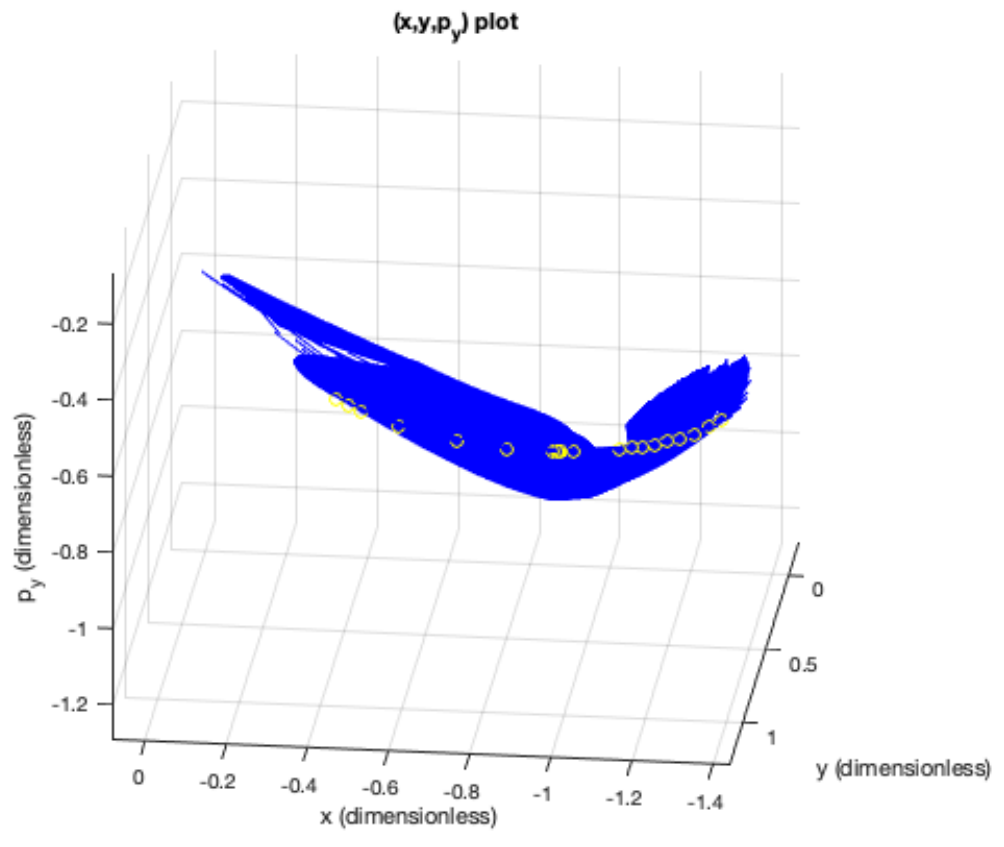}
	          \caption{\label{continuationCartesian} 3D $(x,y,p_{x})$ (left) and $(x,y,p_{y})$ (right) projections of continued intersection points (yellow circles) of various Jupiter-Europa PERTBP 3:4 $W^{u}$ (not shown) with $\omega_{s} = 1.030011$ 5:6 $W^{s}$ (shown in blue)}
	 \end{figure}

\begin{figure}[h]
	\begin{centering}
          \includegraphics[width=0.495\columnwidth]{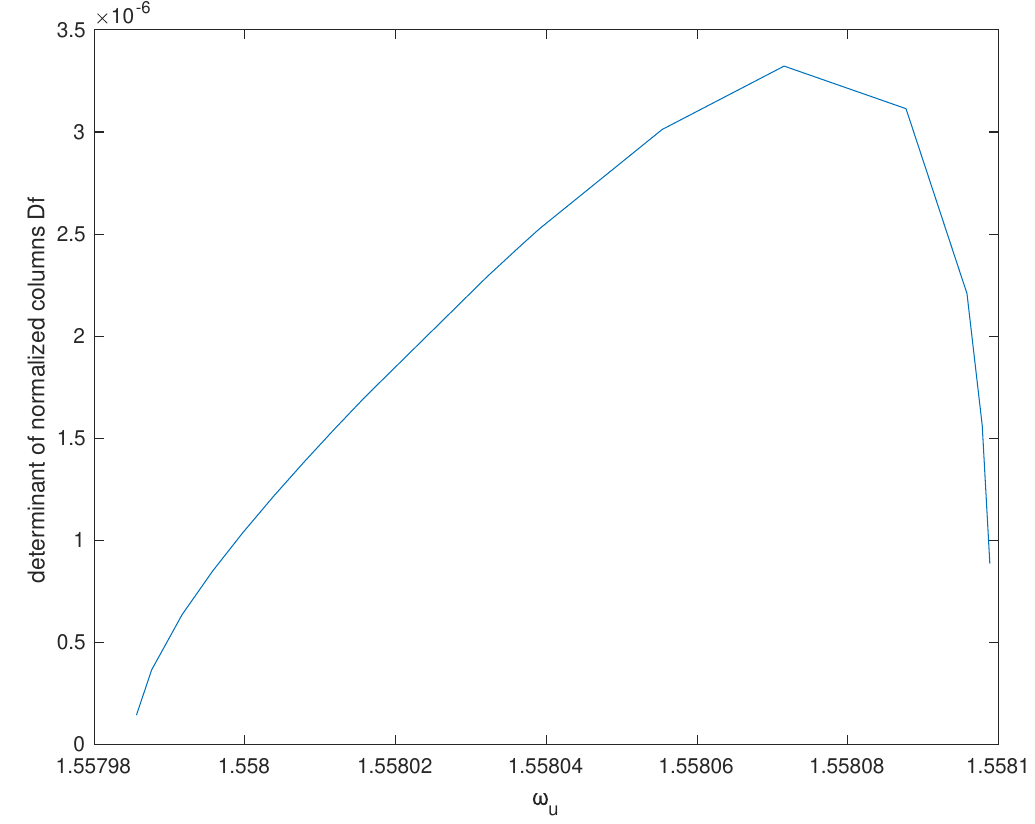}
	          \caption{\label{detVsOmega} Normalized columns $\det Df$ vs $\omega_{u}$}
	  \end{centering}
	 \end{figure}

It is evident that the curves in Figures \ref{solnsThetaVsOmega} and \ref{solnsSvsOmega} are approaching singularities as we approach the lower and upper limits of the range of $\omega_{u}$ values for which we continued the heteroclinic. To see why this is happening, it is instructive to plot the determinant of the $4 \times 4$ matrix derivative, with all columns normalized to unit length, of the function $f$ defined in Equation \eqref{tosolve}. This is shown in Figure \ref{detVsOmega}, where it is immediately clear that this determinant is going to zero as one approaches the upper and lower $\omega_{u}$ limits. From a geometric point of view, this suggests that as $\omega_{u}$ approaches these limits, the intersection of the two manifolds will approach a tangency (non-transversal intersection), past which the manifolds will become locally non-intersecting. From an analytic point of view, this zero determinant will prevent application of the implicit function theorem which is necessary to justify the continuation, which is consistent with the manifolds no longer intersecting past those limits on $\omega_{u}$. 

As a final note, it is worth discussing the significance of the range of 3:4 Jupiter-Europa PERTBP torus rotation numbers $\omega_{u} \in [1.55799, 1.55810]$ from where the heteroclinic connections to the $\omega_{s} = 1.030011$ 5:6 torus were found. Each of these rotation numbers corresponds to a different periodic orbit from the Jupiter-Europa PCRTBP, with the rotation number being determined using Equation \eqref{invariance_fK2} using the fixed perturbation frequency $\Omega_{p}$ and the frequency $\Omega_{1}$ (and thus period) of the periodic orbit. Thus, a range of rotation numbers corresponds to a range of PCRTBP orbit periods. In our case, the 3:4 tori with $\omega_{u} \in [1.55799, 1.55810]$ are continuations of PCRTBP 3:4 periodic orbits with periods $T_{3:4} \in [25.3376, 25.3394]$, while the 5:6 torus with $\omega_{s} = 1.030011$ continues from the PCRTBP orbit of period 38.3281. In the Jupiter-Europa PCRTBP, these 3:4 orbits have Jacobi constants between 3.00234 and 3.00246, while that 5:6 orbit has Jacobi constant 3.0024. Thus, the effect of the PERTBP perturbation on such heteroclinic connections is to allow for (in this case small) changes of energy during a heteroclinic trajectory. This is reminiscent of the phenomenon of Arnol'd diffusion which has been shown\cite{capinski2016} to exist in the PERTBP, where repeated homoclinic connections can be used to effect a significant change in the energy of the spacecraft over time.

\section{Conclusions}

In this paper, we presented a suite of concepts, methods, and tools for finding heteroclinic connections between unstable invariant tori in periodically-perturbed PCRTBP models. Using the idea of layers, we can restrict the connection search to only certain subsets of the manifolds. By generating a discrete mesh of points for the manifolds during globalization, one can combine methods from computer graphics collision-detection algorithms with the massively parallel computing power of modern GPUs for the purpose of rapidly detecting and computing intersections of the meshes. Finally, we showed how to refine the solutions found from the mesh search for greater accuracy, using the manifold parameterizations to enable the application of a differential correction algorithm, which can then be used for numerical continuation through families of tori as well. This produces a range of potential transfer options for space mission trajectory design applications. 

Testing our GPU assisted mesh intersection tools, we saw speedups by a factor of 5-7 as compared to CPU-only tools when using more modern GPU hardware. Even a consumer-grade laptop with a discrete GPU allows for the checking of 14 layers of two manifolds for intersection in just a matter of seconds. This excellent performance, despite the higher system dimension, allows for the feasible investigation of heteroclinic and homoclinic phenomena in such periodically-forced PCRTBP models, where these phenomena can lead to more complex and interesting dynamical behaviors as compared to the PCRTBP (e.g. Arnol'd diffusion \cite{capinski2016}). Our tools thus open up new possibilities for the numerical study of such topics, as well as for finding propellant-free spacecraft trajectories between resonances and other dynamical structures in higher-accuracy dynamical models than the PCRTBP. Trajectories found in such higher-accuracy, periodically-forced PCRTBP models are expected to be easier to numerically continue to the full-accurate ephemeris N-body model used for flying real spacecraft.

Our algorithms are suitable for exploring many different periodically-perturbed PCRTBP models other than the PERTBP as well, such as restricted 4-body models\cite{kumar2023}. The methods could also be adapted to other types of stable/unstable manifolds in such models, such as those of periodic orbits having two stable and two unstable directions (and thus stable and unstable manifolds diffeomorphic to $\mathbb{R}^{2}$ under the stroboscopic map). They could also be extended to the spatial CRTBP in combination with restricting attention to a Poincar\'e section inside a fixed energy surface. Generalizations to higher dimensional tori, manifolds, and systems should also be possible. Another further step in the program would be to develop algorithms chaining the computed heteroclinic connections together. Significant mathematical theory\cite{GideaLS20, FontichM00} exists showing that such chaining of heteroclinics is possible, so we are optimistic that a practical implementation is within reach. Even further steps would be to allow using small amounts of fuel and optimizing among different goals as part of this process. Indeed, there are many possibilities for further development and future work which this methodology presents. 

\section{Acknowledgments}

This work was supported by a NASA Space Technology Research Fellowship under grant no. 80NSSC18K1143. Part of the writing of this article was supported by the National Science Foundation under award no. DMS-2202994. This research was carried out in part at the Jet Propulsion Laboratory, California Institute of Technology, under a contract with the National Aeronautics and Space Administration (80NM0018D0004). The High Performance Computing resources used in this investigation were provided by funding from the JPL Information and Technology Solutions Directorate. R.L.A. was supported through funding by the Multimission Ground System and Services Office (MGSS) in support of the development of the Advanced Multi-Mission Operations System (AMMOS). R.L was supported in part by NSF grant DMS 1800241. Special thanks to Prof Jay Mireles James for suggesting the idea of layers at MSRI.

\bibliographystyle{siamplain}   
\bibliography{references}   

\end{document}